\def \version {2025--01--15}

\documentclass[12pt]{article}

\def \bvt {\begin{verbatim}\footnotesize }

\usepackage{amsmath}
\usepackage{amssymb}
\usepackage{amsthm}
\usepackage{MnSymbol,enumitem}
\usepackage{lineno}
\usepackage{latexsym}

\newcommand{\nev}[1]{{\bf\itshape (#1.)}}
\newtheorem{thm}{Theorem}[section]
\def \btm {\begin{thm}}
\def \etm {\end{thm}}
\newtheorem{prp}[thm]{Proposition}
\def \bpn {\begin{prp}}
\def \epn {\end{prp}}
\newtheorem{lem}[thm]{Lemma}
\def \blm {\begin{lem}}
\def \elm {\end{lem}}
\newtheorem{obs}[thm]{Observation}
\def \bob {\begin{obs}}
\def \eob {\end{obs}}
\newtheorem{rmk}[thm]{Remark}
\def \brm {\begin{rmk}}
\def \erm {\end{rmk}}
\newtheorem{cor}[thm]{Corollary}
\def \bcr {\begin{cor}}
\def \ecr {\end{cor}}
\newtheorem{con}[thm]{Conjecture}
\def \bcj {\begin{con}}
\def \ecj {\end{con}}
\newtheorem{prm}[thm]{Problem}
\def \bpm {\begin{prm}}
\def \epm {\end{prm}}
\newtheorem{dfn}[thm]{Definition}
\def \bdf {\begin{dfn}}
\def \edf {\end{dfn}}
\newtheorem{exa}[thm]{Example}
\def \bex {\begin{exa}}
\def \eex {\end{exa}}

\def \bpf {\begin{proof}}
\def \epf {\end{proof}}

\def \bsk {\bigskip}
\def \msk {\medskip}
\def \ssk {\smallskip}
\def \nin {\noindent}
\def \smin {\setminus}
\def \es {\varnothing}

\newcommand{\floor}[1]{\left\lfloor #1 \right\rfloor}
\newcommand{\ceil}[1]{\left\lceil #1 \right\rceil}
\def \cF {\mathcal{F}}
\def \cfs {\mathcal{F}^*}
\def \cH {\mathcal{H}}
\def \cK {\mathcal{K}}
\def \cM {\mathcal{M}}

\def \fnf {f(n,G|\cF)}
\def \gnf {g(n,G|\cF)}
\def \hnf {g(n,H|\cF)}

\newcommand{\upp}[1]{\langle #1 \rangle}
\newcommand{\nem}[1]{\mathfrak{F}(\neg #1)}
\newcommand{\fs}[1]{\cF\langle #1 \rangle }
\def \fsq {\cF_{\mathrm{sq}}}

\def \fek {\cF_{1,k}}
\def \trif {{\cF^*}}
\def \nti {\to\infty}

\def \knkk {g(n,K_4|\upp{2K_2})}
\def \knkh {g(n,K_4|\upp{K_3})}

\def \arr {anti-Ramsey}
\def \Arr {Anti-Ramsey}

\def \chiF {\chi_{_\cF}}
\def \AP {A}
\def \AR {\mathrm{Ar}}
\def \lar {\mathrm{Lr}}

\newcommand{\Lar}[1]{\lar(n,#1)}

\def \ex {\mathrm{ex}}

\def \Lex {{\sf LEX}}
\def \rtem {{\sf RTEM}}
\def \rtdm {{\sf RTDM}}
\def \rtel {{\sf RTEL}}
\def \rtdl {{\sf RTDL}}

\def \erd {Erd\H os}

\def \Sim {Simonovits}
\def \sos {S\'os}
\def \tur {Tur\'an}

\begin{document}

\title{
 Monochromatic graph decompositions and monochromatic piercing \\
  inspired by anti-Ramsey colorings}
\author{Yair Caro\,\thanks{~Department of Mathematics, University of Haifa-Oranim, Tivon 36006, Israel}
 \and Zsolt Tuza\,\thanks{~HUN-REN Alfr\'ed R\'enyi Institute of Mathematics, Budapest, Hungary} $^,$\thanks{~Department of Computer Science and Systems Technology, University
of Pannonia, Veszpr\'em, Hungary}}
\date{\small Latest update on \version}
\maketitle




\begin{abstract}
\Arr\ theory was initated in 1975 by \erd, Simonovits and \sos,
 inspiring hundreds of publications since then.
Given a graph $G$, the problem is to
 determine the smallest integer $k$
 such that every edge coloring of the complete graph $K_n$
 of order $n$ with at least $k$ colors forces a copy of $G$
 with all edges of distinct colors.
This minimum is denoted by $\AR(n,G)$.

The current paper is the third item of a trilogy in which we
 consider a far-reaching generalization of the \arr\ paradigm.
Two kinds of functions are studied.
For any graph $G$ and family $\cF$ of graphs, introduce:
\begin{itemize}
 \item If $K_2\in \cF$, let $f(n,G|\cF)$ be the smallest integer
  $k$ such that every edge  coloring of $K_n$ with at least
  $k$ colors forces a copy of $G$ in which all color classes
  are members of $\cF$.
 \item If $K_2\notin \cF$, let $g(n,G|\cF)$ be the largest
  integer $k$ for which there exists an edge coloring of $K_n$
  using exactly $k$ colors, such that every copy of $G$ contains
  an induced color class which is a member of $\cF$.
\end{itemize}
Clearly the case $\cF = \{ K_2 \}$ is the classical \arr\
 (aka rainbow) paradigm.

In Part 1 of the trilogy we concentrated mainly on hereditary
 families and developed some of the necessary tools to allow
  lower and upper bounds in this more general setting.
In many cases sharp asymptotic results were obtained.

In Part 2 of the trilogy we concentrated on families
 (not necessarily hereditary) related to popular edge-coloring
 types with parity constraints.

Here in Part 3 we focus on the most general situation, regardless
 of any hereditary consideration or parity constraints.

We show the close relation between the functions $f$ and $g$, and
 develop further useful tools.
We supply dozens of examples (mostly with asymptotically sharp
 bounds) concerning this panoramic view, as well as raise
 many open problems for future research.

To get a glimpse of the spirit and varieties of results proved
 here, we give two illustrative examples of the
 sharp asymptotics obtained.

\msk

\nin
(1)\quad Let $\fsq$ denote the family of all graphs having a
 square number of edges, namely
  $\fsq = \{ G \mid  e(G) = k^2 , \ k \geq 1 \}$.
Then, for every $p \geq 5$, $f(n,K_p|\fsq)  =
 ( 1+o(1))\,\ex(n,K_{\ceil{p/2}})$.

\msk

\nin
(2)\quad Given a graph $G$, denote by $\upp{G}$ the class of
 all graphs containing $G$ as a subgraph.
Assume $p \geq 2q -1$, and let $H$ be any connected graph
 of order $q$.
Then $\ex(n,K_{\ceil{p/(q-1)}}) + \ceil{p/(q-1)} - 1 \leq 
 g(n,K_p|\upp{K_q}) \leq g(n, K_p|\upp{H}) \leq
  (1+o(1))\,\ex(n,K_{\ceil{p/(q-1)}})$.

\bsk

\nin

\bsk

\nin

\end{abstract}

{
\tableofcontents}


\section{Introduction}
\label{s:intro}

\subsection{Anti-Ramsey paradigm and recent developments --- a short overview}
\label{ss:history}

\Arr\ theory was initated in 1975 by \erd, Simonovits and \sos\
 \cite{ESS-75}, inspiring hundreds of publications since then.
The main problem is:
Given a graph $G$, determine the smallest integer $\AR(n,G)$
 such that every edge coloring of the complete graph $K_n$ with
 at least that many colors forces a copy of $G$ with all edges
 getting distinct colors.
The main result of \cite{ESS-75} is:
Let $\chi_e(G) = \min \{ \chi(G - e) : e \in E(G) \}$.
If  $\chi_e(G) = k \geq 3$ then
 $\AR(n,G) = (1+o(1))\,\ex(n,K_k)$, where $\ex(n,G)$ is
 the famous \tur\ number.
The proof relies heavily on the \erd--Stone--Simonovits theorem
 \cite{ES-46,ES-66} that if $\chi(G) = k \geq 3$ then
  $\ex(n,G) = (1+o(1))\,\ex(n,K_k)$. 
Decades later it was finally proved that, for all $n \geq k \geq
 3$,  $\AR(n, K_{k+1}) = \ex(n, K_k) +2$ \cite{MBNL-02,S-04}.
The subject is still very active today, an extensive survey with
 lots of results and references is given in \cite{dyn-surv}.
For the various directions in \arr\ theory we explicitly mention
 here the papers
\cite{{AC-24,AI-08,BGR-15,CLT-09,GR-16,J-02,LLL-24+,MB-06,MBNL-05,Y-21+}};
 some of them are old and some are recent, the list is by no means complete.

Recently the present authors initiated the study of the following more general ``rainbow problem'' \cite{ar-G}.
Let $G$ be a graph, and $\cF$ any given family of graphs.
For every integer $n \geq |G|$, let $f(n, G|\cF)$ denote the smallest integer $k$ such that any edge coloring of $K_n$ with at least $k$ colors forces a copy of $G$ in which each color class induces a member of $\cF$.
Observe that in anti-Ramsey problems each color class is a single edge; i.e., $\cF = \{K_2\}$.

A major result in \cite{ar-G} is:
Let $\cF$ be a hereditary family and let
 $$
   \chiF(G) := \min \{\chi(G-D): D \subset G, \, D \in \cF)\}\,.
 $$
If $\chiF(G) = k \geq 3$ then $f(n,G|\cF) = (1+o(1))\,\ex(n,K_k )$.
Again, a heavy use of the \erd--Stone--Simonovits theorem is inevitable, together with new tools of independent interest, e.g.\ the Independent Transversal Lemma in directed graphs of bounded outdegree.

Also, it has been proved that if $G$ is stable with respect to
 $\cF$, which means $\chiF(G) = \chi(G)$, and if
  $\chi(G) = k \geq 3$ also holds, then
  (regardless of whether or not $\cF$ is hereditary),
 $f(n,G|\cF) = (1+o(1))\,\ex(n,K_k )$.
Many examples of interesting and natural hereditary families and
 the implied results concerning $f(n, G|\cF)$ or $f(n,K_p |\cF)$
 are given in \cite{ar-G}.

Already in \cite{ar-G} we mentioned that certain lower bounds obtained, as well as the case where $\chiF(G) = \chi(G)$, remain true even if $\cF$ is not hereditary but contains $K_2$.

\subsection{Generalized \arr\ monochromatic decomposition and piercing}
\label{ss:4f}

The main purpose of this paper is to consider the following:

Assume that a family $\cF$ is given.
For every graph $G$, we define two functions on the domain of
 all integers $n\geq |G|$.
The first one, $f$, was already the subject of the first \cite{ar-G}
 and second \cite{ar-M} part of the current ``trilogy''.
The other function, $g$, termed \emph{piercing}, will be studied
 in Section \ref{s:pierc} of the present paper.
A main distinction is made depending on whether
 the single edge $K_2$ is, or is not a member of $\cF$.

\msk

\nin
\underline{$\cF$ contains $K_2$\,:}

$\bullet$
Define $\fnf$ as the smallest integer $m$ such that in
 \emph{every} edge coloring of $K_n$ using at least $m$ colors \emph{there exists}
 a copy of $G$ in which \emph{all} color classes are members of $\cF$.

This $f$ is well defined, satisfying $\fnf\leq \binom{n}{2}$,
 as the conditions hold for any $G$ in the rainbow $K_n$.

\msk

\nin
\underline{$\cF$ does not contain $K_2$\,:}

$\bullet$
Define $\gnf$ as the largest integer $m$ such that
 \emph{there exists} a coloring of $K_n$ with $m$ colors, where
 \emph{every} copy of $G$ contains a color class
 which is a member of $\cF$.
If no such coloring exists, then $\gnf:=0$.

A color class in which a subgraph of the considered copy of $G$
 occurs, will be called a \emph{piercing color class}.
For distinct copies of $G$, the piercing color classes are
 allowed to be distinct.
If $\cF$ consists of just one graph, we call this case
 \emph{exact piercing}.

\msk


In the next subsection we will give a more detailed description
 of the structure of the paper.
Briefly saying here, our main goal is to develop the general
 approach to \arr\ theory via detailed analysis of the functions
  $f(n,G|\cF)$ and $g(n,G|\cF)$.
For this purpose we introduce new tools, also
 using classical theorems from extremal graph theory and
 incorporating results about packing graphs.
Doing so, we mostly study the
 case $G = K_p$ with respect to numerous families $\cF$.
The power of our approach is demonstrated
 by a wide range of applications.

\subsection{Structure of paper and sample results}
\label{ss:summary}

Here we sketch the contents of the subsequent sections,
 and mention some representative results or their
 main message, not always quoting all their conditions in
 full detail.
For their strongest forms please see the corresponding theorems.

Section \ref{s:col-typ} describes some basic inequalities and
 coloring patterns that are used along the paper.
Among the simple relations we mention the following important
 principle of complementarity:

$\bullet$\quad
If $K_2 \in \cF$ and $G \notin \cF$, then
 $g(n,G|\overline{\cF}) = f(n,G|\cF) - 1$.
  (Proposition~\ref{p:compl})

This allows deriving results on one of
 $f$ and $g$, and obtaining an immediate result on the other function.
It is highly beneficial, as the methods developed for estimating
 the two functions are in general much different from each other.
We apply this relation several times along Section~\ref{s:pierc},
 and a general consequence (Proposition~\ref{p:compl-c,d}) is
  also derived from it in the concluding section.
However, the reader may discover that the complementarity
 principle can be applied in many more situations where we
 do not mention it one by one along the text.

In Section \ref{ss:quad-UB} we collect some important results
 by which efficient upper bounds on $f(n,G|\cF)$ can be obtained. 
Among them are the \erd--Rado theorem, the \erd--Stone--Simonovits
 theorem, the Rainbow Cut Lemma that we developed in \cite{ar-G},
  and the Canonical Cut Lemma proved here as Theorem \ref{t:CCCL},
 that will be a major tool in the present paper. 

Section \ref{s:quadra} concerns numerous families
 for which $f(n,K_p|\cF)$ has quadratic growth order.
Among those $\cF$ are the following families:
 all regular graphs;
 all graphs except the complete graphs 
$K_r$ with $r\geq q$ for any fixed $q\geq 3$ if $p \geq 2q-1$;
 all graphs except the starts $K_{1,r}$ with $r\geq q$ for any fixed $q\geq 2$ if $p \geq 2q+1$;
 $K_2$ together with all graphs containing at least one cycle;
 all graphs of maximum degree at most $\Delta$ for any fixed $\Delta$;
 families with connectivity constraints; and more.
From the many tight asymptotic results we mention two examples: 

$\bullet$\quad
If $p$ and $r$ are integers with either $p\geq 4$ and $r=1$
 or $p\geq 2r$ and $r\geq 2$, and $\cF=\{mK_2 \mid 1\leq m\leq r\}$,
 then $f(n,K_p|\cF) = (1+o(1))\,\ex(n,K_{p-r})$.
The asymptotic formula remains the same if $\cF$ is supplemented
 with all triangle-free graphs. (Theorem \ref{t:r-match})

$\bullet$\quad
If $\cF$ is the family of all graphs on at most $k$ vertices
 ($k\geq 2$), and
 $p\geq k+2$, then $f(n,K_p|\cF) = (1+o(1))\,\ex(n,K_{p-k+1})$.
  (Theorem \ref{t:k-small})

Section \ref{ss:stab} concerns a certain kind of stability.
As a very general result beyond $K_p$,
  we prove a sufficient condition for the
 existence of graphs which are chromatically stable with respect
 to a family $\cF$ and make $f$ grow proportionally to $n^2$
 (more precisely $f(n,G|\cF)\sim \ex(n,K_r)$ for a given $r\geq 3$).

$\bullet$\quad
If the family $\cF$ contains $K_2$, and
 $\max\{ e(D) : D\in \cF, \, |V(D)|=n\} < n^2/k^2$ holds
 for every $n$ where $k\geq 3$ is a fixed integer, then
 there exist infinitely many graphs $G$ such that
 $\chiF(G)=\chi(G)=k$ and all those $G$ have
 $f(n,G|\cF) = (1+o(1))\,\ex(n,K_{k})$.  (Theorem \ref{t:k2n2})

Section \ref{s:excl1}  concerns the behaviour of $f(n,K_p|\cF)$
 where $\cF$ is obtained from the class of all graphs 
 by omitting just one graph $H$.
These densest nontrivial families will be major tools to handle
 instances of exact piercing in Section~\ref{s:pierc}.
There is a surprising difference if the deleted graph $H$
 is a star or any other graph:

$\bullet$\quad
If $\cF$ is obtained by excluding the star $K_{1,t}$, then
 $f(n,K_p|\cF) = (1+o(1))\,\ex(n,K_{\ceil{p/t}})$ if $p\geq 2t+1$,
 and $o(n^2)$ if $t+1 \leq p \leq 2t$.  (Theorem \ref{t:excl-star})

$\bullet$\quad
If $\cF$ is obtained by excluding a non-star $G\neq K_p$ of
 order at most $p$, then $f(n,K_p|\cF) = 1$ for all
 $n > n_0(p)$.  (Theorem~\ref{t:excl-nonstar})

Section \ref{s:quad-spec} concerns the asymptotic sharp
 determination of $f(n,K_p|\cF)$ where the members of $\cF$ are
 the graphs whose number of edges belongs to a prescribed
 integer sequence $S$.
Among the considered cases are various arithmetic progressions,
 and also the much sparser sequence of all squares.

$\bullet$\quad
If $p\geq 5$ and $\cF$ is generated by $S=\{ k^2 \mid k\geq 1\}$,
 or $p\geq 4k-1$ and $S=\{m \mid m\equiv 1$ (mod $k)\}$, then
 $f(n,G|\cF) = (1+o(1))\,\ex(n,K_{\ceil{p/2}})$.
  (Theorems \ref{t:square} and \ref{t:1-mod-k})

Section \ref{s:subqu} concerns the existence of families $\cF$
 for which $f(n,G|\cF) \leq h(n)$ where $h(n)$ is a
 prescribed function that ranges from a constant function to
 a quadratic growth function.

$\bullet$\quad
If $\cF$ is obtained by excluding all graphs $H$ with
 $\AR(n,H)\leq h(n)$ where $h(n)$ is a positive integer function,
 then $f(n,G|\cF) \leq h(n)$ if $G\notin \cF$, and $f(n,G|\cF)=1$
 for all $n>n_0(G)$ if $G\in \cF$.   (Proposition \ref{p:small-AR})

Section \ref{s:pierc} constitutes a major contribution of this
 paper concerning the notion of piercing and the function $g(n,G|\cF)$.
Before entering into details, let us emphasize that the
 complementary relation between $f$ and $g$ allows to derive 
 asymptotic results for numerous families in either direction. 

After a short part in which we prove subquadratic growth under
 certain conditions, a  subsection is devoted to $g(n,K_p|\cF)$
 where $\cF=\upp{H}$ is the family of all graphs containing a
 specified small graph $H$.
This track is done in parallel to \arr\ type results given in
 \cite{BGR-15} and \cite{GR-16}.
Its aim is to show also the extra efforts needed to obtain exact
 formulas, rather than asymptotic sharp bounds.

$\bullet$\quad
We have $g(n,K_4-e|\upp{K_3})=1$ and $g(n,K_4|\upp{K_3})=n$ for
 all $n$,
 whereas $g(n,K_4|\upp{P_3})=\Theta(n^{3/2})$ and
 $g(n,K_4|\upp{\{K_{1,3},P_4\}})=\Theta(n^{3/2})$ for $n\nti$.
 (Theorems \ref{t:uppK3} and \ref{t:uppP3})

We then develop the method of graph packing and anti-packing as a
 major tool to derive bounds on $g(n,G|\cF)$. 
This approach allows us to prove numerous asymptotic
 sharp results.
We quote the following ones here:

$\bullet$\quad
If a graph $G$ of order $p$ has no isolated vertices, then
 $g(n,K_p|\upp{G})=O(n)$; and if it has isolates but its
 maximum degree is at least $p/2$, then $g(n,K_p|\upp{G})=
 O(n^{2-\epsilon})$ for some $\epsilon=\epsilon(G)>0$ as
 $n\nti$. (Theorem \ref{t:minmaxdeg})

$\bullet$\quad
For every $p\geq 5$, $q\geq 3$, and $n>n_0(p)$:
 $g(n,K_p|\upp{K_q}) = (p-q)n - \binom{p-q+1}{2} + 1$ if $q < p\leq 2q-3$, an exact formula for large $n$;
 $g(n,K_p|\upp{K_q}) = (p-q)n - O((p-q)^2)$ if $p=2q-2$; and
 $g(n,K_p|\upp{K_q}) = (1+o(1))\,\ex(n,K_{\ceil{p/(q-1)}})$
  if $p\geq 2q-1$. (Theorems 
  \ref{t:gnpq}, \ref{t:kpq-quad}, and \ref{t:2q-2})

$\bullet$\quad
If $p\geq 2t+3\geq 5$ and $\cF$ is the family of all graphs which
 do not contain $K_{p-t}$, then
 $f(n,K_p|\cF) = t\cdot (n-t) + \binom{t}{2} + 2$.
  (Theorem \ref{t:no-k-p-t})

$\bullet$\quad
If $p\geq 2q-1\geq 5$ and $\cF$ is the family of all $K_q$-free
 graphs, then
 $f(n,K_p|\cF) = (1+o(1))\,\ex(n,K_{\ceil{p/(q-1)}})$.
  (Theorem \ref{t:p>2q})

$\bullet$\quad
If $p\geq 5$ and $n>n_0(p)$, then $g(n,K_p|\upp{C_p}) =
 \ceil{\frac{p-3}{2}\,n} + 1$.  (Theorem~\ref{t:kpcp})

Section \ref{s:conclude} contains a concluding summary of the
 trilogy, i.e., based upon \cite{ar-G}, \cite{ar-M}, and the
  current paper.
The section also provides an extensive list of many problems,
 with the intention to motivate further development of the
  subject in future research.

\subsection{Standard definitions and some special notation}

Most of our notation is standard.
For a graph $G$ we denote by $|G|$ the number of vertices and
 by $e(G)$ the number of edges.
The degree of vertex $v$ is denoted by $d_G(v)$, and is
 abbreviated as $d(v)$ when $G$ is understood.
Also, as usual, $\delta(G)$ and $\Delta(G)$ stand for the
 minimum degree and maximum degree of $G$, respectively.
The path, the cycle, and the complete graph on $k$ vertices is
 denoted by $P_k$, $C_k$, and $K_k$, respectively.
Further, $K_{a,b}$ denotes the complete bipartite graph with
 $a$ and $b$ vertices in its vertex classes.

A standard notation in extremal graph theory is
 $\ex(n,F)$ for the \tur\ number of graph $F$, defined as the
 maximum number of edges in a graph of order $n$ that does not
 contain any subgraph isomorphic to $F$.
In case of $F=K_{m+1}$ the unique extremal graph---the \tur\
 graph---is denoted by $T_{n,m}$.
It is the complete $m$-partite graph in which the sizes of
 any two vertex classes differ by at most 1.

For the sake of easier reading we also collect here some
 non-standard notation that will be important later,
 although they will be defined in the sequel.
We consider the following huge families of graphs:
 \begin{itemize}
  \item $\mathfrak{F}$ --- all
 finite graphs without isolated vertices;
  \item $\nem{F}$ --- omitting the single given graph $F$
   from $\mathfrak{F}$;
  \item $\upp{F}$ --- all graphs containing the given
   graph $F$ as a subgraph;
  \item $\fs{S}$ --- all graphs $H$ with $e(H)\in S$, where $S$
   is a given set of positive integers.
 \end{itemize}

The families $\upp{F}$ of graphs generated by a
 single graph $F$ are of great
 interest in the context of function $g$.
A part of our study we will concentrate on $G=K_p$ and $F=K_q$.
For this purpose we introduce the notation
 $$
   g(n,p,q) := g(n,K_p|\upp{K_q})
 $$
 where $p>q\geq 3$ will naturally be assumed.

Beside the \arr\ function $\AR(n,G)$ we shall also apply its
 local version $\lar(n,G)$ in the sequel.
It is defined as the smallest integer $k$ such that
 every edge coloring of $K_n$ with at least $k$ colors
 contains a properly edge-colored copy of $G$.
This function was introduced and studied in
 \cite{ar-G,ar-M}.

\section{Some basic inequalities and coloring types}
\label{s:col-typ}

We begin with some simple inequalities.

\bpn
\nev{Monotonicity}   \label{p:monoton}
\begin{itemize}
 \item[$(i)$] If\/ $K_2\in \cF_1\subset \cF_2$, then
   for every\/ $G$ we have\/ $f(n,G|\cF_1) \geq f(n,G|\cF_2)$.
 \item[$(ii)$] If\/ $\cF$ is hereditary,
   and\/ $H\subset G$, then\/
  $f(n,H|\cF) \leq f(n,G|\cF)$.
 \item[$(iii)$] If\/ $\cF_1\subset \cF_2$ and\/ $K_2\notin \cF_2$, then
   for every\/ $G$ we have\/ $g(n,G|\cF_1) \leq g(n,G|\cF_2)$.
 \item[$(iv)$] If\/ $H\subset G$ and\/ $Q\subset H$, then\/
  $g(n,H|\upp{Q}) \leq g(n,G|\upp{Q})$.
 More generally, if\/ $H\subset G$, then\/
  $g(n,H|\cF) \leq g(n,G|\cF)$ holds for every\/ $\cF$
  ($K_2\notin\cF$) closed under supergraphs; i.e.,
  if\/ $F\in\cF$ and\/ $F\subset F'$, then also\/ $F'\in\cF$.
 \item[$(v)$] If\/ $Q_1\subset Q_2\subset G$, then\/
  $g(n,G|\upp{Q_1}) \geq g(n,G|\upp{Q_2})$.
\end{itemize}
\epn

\bpf
$(i)$ \ Since every member of $\cF_1$ also belongs to $\cF_2$,
 a copy of $G$ having an $\cF_1$-decomposition by its color
 classes satisfies the same requirement for $\cF_2$ as well.

\ssk

\nin
$(ii)$ \ If $G'\subset K_n$ is a copy of $G$ valid for $\cF$,
 then any subgraph $H'\subset G'$ isomorphic to $H$ inherits a
 coloring in which the monochromatic classes are subsets of
 those in $G'$, hence are members of $\cF$ because
 $\cF$ is hereditary.

\ssk

\nin
$(iii)$ \ Using $g(n,G|\cF_1)$ colors, every copy of $G$ contains
 a color class that belongs to $\cF_1$, and by assumption this
 subgraph is also a member of $\cF_2$.

\ssk

\nin
$(iv)$ \ The inequality trivially holds if $\hnf=0$.
If $\hnf=k>0$, let $\psi$ be an edge $k$-coloring of $K_n$ that
 pierces every copy of $H$.
Any $G'\cong G$ in $K_n$ contains a $H'\cong H$.
Let $c$ be the color that pierces $H'$ under $\psi$ with a
 subgraph $F\in\cF$.
Then the $c$-colored subgraph $F'$ of $G'$ contains $F$, hence
 $F'\in\cF$ and therefore also $G'$ is pierced.

\ssk

\nin
$(iv)$ \ If $\psi$ is an edge coloring of $K_n$ that pierces
 every copy of $G$, for $\upp{Q_2}$, then it is feasible also
 for $\upp{Q_1}$, because any monochromatic supergraph of $Q_2$
 contains $Q_1$ as a subgraph, hence is a member of $\upp{Q_1}$.
So, $\upp{Q_1}$ admits at least as many colors as $\upp{Q_2}$ does.
\epf

There is also a complementary relation between the
 functions $f$ and $g$.
Then, for any family $\cF$ of graphs with $K_2\in\cF$, define its
 complementary family $\overline{\cF} := \mathfrak{F} \smin \cF$,
 obrained from the universe $\mathfrak{F}$ of all finite graphs
 without isolated vertices by removing all members of $\cF$.

\bpn
\nev{Family complementation} \  \label{p:compl}
\begin{itemize}
 \item[$(i)$] For every\/ $\cF$ with\/ $K_2\in\cF$, we have\/
 $g(n,G | \overline{\cF}) < f(n,G |\cF)$.
 \item[$(ii)$] If\/ $K_2\in\cF$ and\/ $G\notin \cF$, then\/
 $g(n,G | \overline{\cF}) = f(n,G |\cF)-1$.
\end{itemize}
\epn

\bpf
Let $\psi$ be any edge coloring of $K_n$ using $k$ colors, where
 $k\geq f(n,G |\cF)$.
By the definition of $f$, there exists a copy $G'$ of $G$
 such that every color class of $G'$ belongs to $\cF$.
Hence, this $G'$ does not contain any color class from
 $\overline{\cF}$, violating the conditions on $g$, therefore
  implying that $g(n,G | \overline{\cF}) = k$ is impossible.
Since this impossibility holds for all $k\geq f(n,G |\cF)$, the
 claimed inequality $(i)$ follows.

For the proof of $(ii)$, let us denote $k:=f(n,G |\cF)$.
By the definition of $f(n,G |\cF)$, there exists an edge
 coloring $\psi$ of $K_n$ using $k-1$ colors, in which
 every copy of $G$ induces a color class not belonging to $\cF$.
It follows that every copy of $G$ contains a color class from
 $\overline{\cF}$, thus $g(n,G | \overline{\cF})\geq k-1$.
Applying $(i)$ we obtain $g(n,G | \overline{\cF})=k-1$,
 hence proving $(ii)$.
\epf

Before proceeding with the next section, the following two simple
 remarks should be taken into account for the rest of the paper.

\brm
An important warning is that if\/ $G\in\cF$---that is, if\/
 $\chiF(G)=1$ holds---one cannot conclude\/ $f(n,G|\cF)=1$.
The fundamental class of counterexamples is of the form\/
 $\cF=\{G,K_2\}$, assuming that\/ $G$ contains a cycle.
Define the edge coloring\/ $\psi$ of\/ $K_n$ as\/ $\psi(v_i,v_j)=j-1$
 for all\/ $1\leq i<j\leq n$ (see the\/ \Lex\ pattern in
 Definition \ref{d:homo} below).
It implies\/ $f(n,G|\cF) \geq n$ because the vertex of highest
 index in a cycle is incident with a star color class, other than\/
 $K_2$, induced by\/ $\psi$ in\/ $G$.
\erm

\brm   \label{r:g-pos}
\
\begin{itemize}
 \item[$(i)$] A necessary condition for\/ $\gnf\geq 1$ is to have a subgraph\/
 $F\subseteq G$ with\/ $F\in\cF$.
 \item[$(ii)$] A sufficient condition for\/ $\gnf\geq 1$ is to have a subgraph\/
 $F\subseteq G$ with\/ $\upp{F}\subseteq\cF$, where\/ $\upp{F}$
  denotes the family of graphs containing\/ $F$ as a subgraph.
\end{itemize}
\erm

\subsection{Some fundamental types of colorings of $K_n$}
\label{ss:patterns}

Typically, lower bounds on the considered functions are obtained
 by explicit constructions, which use standard coloring patterns
 or their combinations on the edges of $K_n$.

\bdf \nev{Homogeneous patterns}   \label{d:homo} \
The following three patterns are the standard homogeneous ones.
\begin{itemize}
 \item Monochromatic: all edges of\/ $K_n$ have the same color.
 \item Rainbow: each edge of\/ $K_n$ has its private color.
 \item \Lex : the vertices are in linear order\/ $v_1,v_2,\dots,v_n$,
  and each edge\/ $v_iv_j$ with\/ $1\leq i<j\leq n$ has color\/ $j-1$.
\end{itemize}
\edf

\paragraph{Compound patterns.}

The general approach to lower-bound constructions with
 fast-growing (superlinear) numbers of colors is to take a
 large rainbow-colored graph $H\subset K_n$ and use one or
 relatively few additional colors in its complement $\overline{H}$.

In cases where the maximum number of colors is $o(n^2)$,
 usually the structure of $H$ is not explicitly known, and
 we take a monochromatic complement.
In situations with quadratic growth, it is convenient to let $H$
 be a \tur\ graph extremal for a certain $K_q$; that is, the
 complete $(q-1)$-partite graph $T_{n,q-1}$ where each of the
 vertex classes $V_1,\dots,V_{q-1}$ has $\floor{\frac{n}{q-1}}$
 or $\ceil{\frac{n}{q-1}}$ vertices.
Then one has the freedom to choose an edge coloring inside each $V_i$.

The original \arr\ problem was solved asymptotically in \erd, \Sim\
 and \sos's paper \cite{ESS-75} by taking monochromatic complement.
In our general theory, however, in many situations it is
 useful---and sometimes substantial, unavoidable---to
 apply \Lex\ inside some or all of the $V_i$.
So, the four general schemes for constructions yielding quadratic growth are:
\begin{itemize}
 \item \rtem$(n,q-1)$\,: rainbow $T_{n,q-1}$ + monochromatic $K_{|V_i|}$
  using the same color inside all $V_i$;
 \item \rtdm$(n,q-1)$\,: rainbow $T_{n,q-1}$ + monochromatic $K_{|V_i|}$
  inside all $V_i$, whose colors are mutually distinct;
 \item \rtel$(n,q-1)$\,: rainbow $T_{n,q-1}$ + \Lex\ inside each $V_i$,
  using the same colors for \Lex\ in every vertex class.
 \item \rtdl$(n,q-1)$\,: rainbow $T_{n,q-1}$ + \Lex\ inside each $V_i$,
  each of those $V_i$ having its private \Lex\ colors.
\end{itemize}
Hence, to make it easy to remember: RT = rainbow \tur,
 M = monochromatic parts, L = \Lex-colored parts,
 E = equally colored, D = distinct colors are used.
We include the number $n$ of vertices in the notation, and the
 second parameter denotes the number of parts (vertex classes)
 of the rainbow \tur\ graph involved.

\section{Tools for asymptotic upper bounds}
\label{ss:quad-UB}

Here we describe a method that will be applied many times in
 proofs of asymptotic upper bounds.
We take the following three major ingredients.

\btm \nev{\erd--Rado Canonical Theorem \cite{ER-50}} \  \label{t:E--R}
For every\/ $t\geq 2$ there exists an integer\/ $n_0(t)$ such
 that, for every\/ $n\geq n_0(t)$, in every edge coloring of\/
 $K_n$, there exists a complete subgraph\/ $K_t$
 which is either monochromatic or rainbow or\/ \Lex-colored.
\etm

\btm \nev{\erd--Stone--\Sim\ Theorem \cite{ES-46,ES-66}} \  \label{t:E-St-S}
For every\/ $a\geq 2$ and every\/ $r\geq 3$, the \tur\ number of
 the complete\/ $r$-partite graph\/ $K_{a,\dots,a}$ satisfies\/
 $\ex(K_{a,\dots,a}) = (1+o(1))\,\ex(n, K_r) = 
 \frac{r-2}{2r-2} n^2 + o(n^2)$ as\/ $n\to\infty$.
\etm

As an immediate consequence, the same asymptotic holds for the \tur\ number
 of any $r$-chromatic graph ($r\geq 3$),
 but we won't need this seemingly more general fact.

\btm \nev{Rainbow Cut Lemma \cite{ar-G}}   \label{t:RCL}
 \
For any integers\/ $b\geq 2$ and\/ $r\geq 2$
 there exists an\/ $a=a(b,r)$ with the following property.
If\/ $\psi$ is an edge coloring of\/ $K_{ar}$,
 and\/ $K \cong K_{a,\dots,a}$ is a rainbow spanning subgraph of\/
 $K_{ar}$ under\/ $\psi$, then $K$ contains a (rainbow)
 subgraph\/ $K' \cong K_{b,\dots,b}$ of order\/ $br$ such that
 the colors of edges inside the vertex classes of\/ $K'$
 in\/ $K_{ar}$ do not occur in\/ $K$.
\etm

Combining the above, the following result can be derived.
For $r\geq 3$ it is a strong tool to prove asymptotically tight
 estimates; and if $r=2$, then the involved \tur\ function
 is zero and the obtained bound is the subquadratic $h(n)$,
 yielding an arbitrarily large ``clean'' complete bipartite
 rainbow subgraph derived from as few as $o(n^2)$ colors.

\btm \nev{Canonical Cut Lemma}   \label{t:CCCL}
 \
For any integers\/ $r\geq 2$ and\/ $t\geq 2$ there exists a
 function\/ $h(n) = h(n,r,t)$ of subquadratic growth rate\/
  $h(n) = o(n^2)$ with the following property.
In every edge coloring\/ $\psi$ of\/ $K_n$ with more than\/
 $\ex(n,K_r) + h(n)$ colors, there occurs a complete\/
 $r$-partite graph\/ $K\cong K_{t,\dots,t}$, say with vertex
  classes\/ $V_1,\dots,V_r$, such that
  \begin{itemize}
   \item[$(i)$] $K$ is rainbow-colored,
   \item[$(ii)$] for\/ $i=1,\dots,r$ each\/ $V_i$ induces a\/ $K_t$
    which is either monochromatic or rainbow or\/ \Lex-colored,
   \item[$(iii)$] no color of\/ $K$ is repeated inside any\/ $V_i$.
  \end{itemize}
\etm

\bpf
Given $t$, we choose $b=n_0(t)$ according to Theorem \ref{t:E--R},
 and then take $a=a(b,r)$ satisfying the requirements of
 Theorem \ref{t:RCL}.
For fixed $r$ and $t$ also $a$ is fixed, hence by Theorem \ref{t:E-St-S}
 the difference $\ex(n,K_{a,\dots,a}) - \ex(n,K_r)$ can be
 expressed as a subquadratic function $h(n,r,t)$.
Then more than $\ex(n,K_r) + h(n,r,t)$ colors always imply the
 presence of a rainbow $K_{a,\dots,a}$;
 it can be cleaned to a rainbow $K' \cong K_{b,\dots,b}$
 with the additional restriction that no color of
 any $e\in E(K')$ appears in $\binom{V(K')}{2} \smin E(K')$;
 and then each vertex class of $K_{b,\dots,b}$ contains a
 complete subgraph of order $t$, which is either monochromatic
 or rainbow or \Lex-colored.
\epf

With reference to part $(iii)$ of this theorem, we will use the following terminology.

\bdf  \nev{Clean multipartite graph} \
A rainbow complete multipartite graph\/ $K\subset K_n$ with
 vertex classes\/ $V_1,\dots,V_r$ is said to be clean if
 none of the\/ $V_i$ induces any edge whose color occurs in\/ $K$.
\edf

In many situations we shall ensure further that the vertex classes
 have homogeneous colorings (monochromatic, rainbow, or \Lex);
 but this will not be a requirement for being ``clean''.

\section{Quadratic growth: The role of exclusions and structural requirements}
\label{s:quadra}

Here we present families $\cF$ for which there exist graphs
 $G$ with $f(n,G|\cF)\geq cn^2$ for some $c>0$.
In most cases we shall achieve much more, namely
 tight asymptotic bounds.
The general principle in obtaining effective lower bounds is
 to apply some of the composite patterns
\begin{center}
 \rtem$(n,m)$,\quad \rtdm$(n,m)$,\quad \rtel$(n,m)$,\quad \rtdl$(n,m)$
\end{center}
 introduced in Section \ref{ss:patterns}.
For the readers' convenience we recall with emphasis that the
 core of these constructions is a rainbow-colored
 complete $m$-partite \tur\ graph $R\cong T_{n,m}$, and none of
  the $\ex(n,K_{m+1})\sim \frac{m-1}{2m}n^2$ colors of $R$ is
  repeated inside any vertex class $V_i$;
 i.e., the color set of $K_n-R$ is disjoint from that of $R$.

The fewest and the most colors occur in \rtem$(n,m)$ and
 \rtdl$(n,m)$, respectively; however, the difference is just
 $n-m-1$, which is negligible in comparison to $\Theta(n^2)$
 as $n$ gets large.
The distinctive properties of the four variants---as expressed
 with the suffixes {\sf EM}, {\sf DM}, {\sf EL}, {\sf DL}---become
 relevant when particular families $\cF$ are considered.

\subsection{Three case studies}

In this subsection we collect three simple and very natural graph classes.

\paragraph{(a) \ $\cF$ = \{regular graphs\}.}

Clearly, $K_5$ can be decomposed into two copies of $C_5$;
 but it is regular, hence $\chi(K_5|\cF) = 1$. 
However, $f(n, K_5|\cF)$ is a quadratic function of $n$,
 as demonstrated by the \rtdl$(n,2)$ coloring.
Then, in one side there must be at least three vertices of $K_5$
 that form $K_3$ (or $K_4$, or the entire $K_5$), colored by
 \Lex\ which forbids monochromatic regular graphs except $K_2$.  Hence $f(n,K_5|\cF) \geq \ex(n,K_3) + n -2
 = \left\lfloor \frac{1}{4}(n^2+4n-8) \right\rfloor$.  

On the other hand, $\lar(n,K_5) = (1+o(1))\,\ex(n,K_3)$ where
 the function $\lar$ requires proper edge coloring;
 i.e., all monochromatic edge classes are
 1-regular, so that $f(n,K_5|\cF) \leq \lar(n,K_5)$ holds.
As a consequence, the asymptotics of
 $f(n,K_5|\cF) = (1+o(1))\,\ex(n,K_3)$ is determined.

Exactly the same proof gives $f(n, K_ p |\cF) >
 \ex(n, K_{\lceil p/2 \rceil})$ using either of the
 \rtel$(n,\lceil p/2 \rceil-1)$ or
  \rtdl$(n,\lceil p/2 \rceil-1)$ coloring pattern.
Then Theorem 5.2 in \cite{ar-G}, stating
 $\lar(n,K_p)=(1+o(1))\,\ex(n, K_{\lceil p/2 \rceil})$,
 gives a matching asymptotic upper bound.
Clearly, although $G \in\cF$, the fact $\chiF(G) = 1$ is useless,
 rather we applied the property of \Lex\ that it allows only stars as monochromatic classes; and for $K_p$ we computed the order $\lceil p/2 \rceil$ that can make use of this property.   

We will show that this can be done in a wider setting, where,
 given $\cF$, a certain large enough induced subgraph of $G$ 
 is anti-\Lex\ with regard to $\cF$. 

\bsk

Before some steps towards the solution of this problem, let us
 give two further case studies as negative examples;
 i.e., graph classes that do admit quadratic growth.

\paragraph{(b) \ Omitting large complete graphs.}

Let $t\geq 3$ be a fixed integer, and let
 $\cF$ be obtained from the family of all graphs by
 omitting just the complete graphs $K_r$ for all $r\geq t$.
We now take the \rtdm$(n,2)$ coloring of $K_n$.
This yields a quadratic lower bound for every complete
 graph $K_p$ of order $p\geq 2t-1$.
Indeed, since $p\geq 2t-1$, any copy of $K_p$ has at least $t$
 vertices in one of the monochromatic classes (and its color is
 distinct from all the other color classes), hence forces a
 color class $K_r$, $r \geq t$, which is not a member of $\cF$.

\paragraph{(c) \ Omitting large stars.}

Let $t\geq 2$ be a fixed integer, and let
 $\cF$ be obtained from the family of all graphs by
 omitting just the stars $K_{1,r}$ for all $r\geq t$.
Take now the edge coloring \rtdl$(n,2)$ of $K_n$.
This coloring yields a quadratic lower bound for every complete
 graph $K_p$ of order $p\geq 2t+1$.
Indeed, any copy of $K_p$ has at least $t+1$ vertices in
 one of the \Lex-colored classes, and its edges
 are split there into monochromatic stars whose colors are
 distinct from all the other color classes.
The maximum degree among those stars is at least
 $\ceil{p/2} - 1 \geq t$, hence at vertex $t+1$ in the \Lex\
  order there is a color class $K_{1,t}$ which is not in $\cF$.

\subsection{Cycles}
\label{ss:quad-LB}

Let $\cF$ be any family of graphs each containing a cycle.
As already explained above, in order to allow decomposition for
 all graphs, we set $\cfs := \cF \cup \{K_2\}$.

\brm
If\/ $G$ is acyclic, then the only decomposition into members
 of\/ $\cF^*$ is to make\/ $G$ rainbow (i.e., all color classes are\/
 $K_2$), and then\/ $f(n,G|\cF^*) = \AR(n,G)$.
\erm

For this reason we will assume here that $G$ is not acyclic.
Note that the addition of $K_2$ is necessary even if $G$ is
 decomposable into members of $\cF$, because a rainbow $K_n$
 does not admit any decomposition other than with classes $K_2$.

It follows from the above discussion that, for $p \geq 5$ and
 $G=K_p$, the \rtel$(n,2)$ or \rtdl$(n,2)$ coloring of $K_n$
 with its more than $\ex(n,K_3)$ colors
 forces $f(n,G|\cfs)$ to be quadratic, because any $K_r$ with
 $r \geq 3$ in \Lex\ induces a color class $K_{1,r-1}$,
 which is not a member of the current $\cfs$.

This example can be raised to a more general level.
For any non-forest graph $G$ let us define $\mu(G)$ as the
 smallest integer $k$ such that every induced subgraph of $G$
 on $k$ vertices contains a cycle; and if $G$ is acyclic,
 then let $\mu(G) = \infty$.

\blm
If\/ $G$ contains a cycle, then\/ $\mu(G) \leq 2\alpha(G) +1$.
\elm

\bpf
Suppose there is an induced subgraph $H$ of $G$ on
 $2\alpha(G) +1$ vertices containing no cycle.
Then $H$ is acyclic, hence bipartite, and $\alpha(G) \geq
 \lceil \, |H|/2 \rceil = \alpha(G) +1$, a contradiction.
\epf

\bpn
Suppose\/ $m \geq 2$ and\/ $\lceil \, |G|/m \rceil \geq \mu(G)$.
Then\/ $f(n,G|\cfs) \geq \ex(n,K_{m+1}) + n-m+1$.
\epn

\bpf
Consider the \rtdl$(n,m)$ coloring built upon the rainbow \tur\
 graph $T_{n,m}$ with $\ex(n,K_{m+1})$ edges,
 using $n-m$ additional colors.
In any copy $G'$ of $G$, there is a class of $T_{n,m}$ containing
 a subgraph $H\subset G'$ with at least $\mu(G)$ vertices.
By definition, this $H$ contains a cycle, say $C$,
 which is \Lex-colored allowing only stars as color classes.
Consider the vertex $v\in V(C)$ of highest index under \Lex.
The edges incident with $v$ in $H$ form a monochromatic star
 other than $K_2$ (as it contains at least two edges of $C$),
 and no other edge of $H$---neither of $G$---has this color.
Thus, the subgraph induced by this color class in $G$ is
 not a member of $\cfs$, and consequently no copy of $G$ is
 $\cfs$-colored.
\epf

We have to note that the role of cycles is substantial.
This fact is demonstrated also by the next observation, which
 was first proved in \cite{AI-08}, using induction.
We include here a short argument that applies the well known
 BFS algorithm.

\bpn   \label{p:tree-lex}
Every forest\/ $T$
 can be packed in\/ \Lex\ so that\/ $T$ is rainbow.
\epn

\bpf
We may assume that $T$ is connected (otherwise extending $T$
 to a spanning tree does not make the task of rainbow
 embedding easier).
Fix a root of $T$ and place it as the lowest vertex of the
 \Lex\ order.
Then pack the vertices step by step upward using Breadth First
 Search.
Then at each step every new vertex has only one lower neighbor,
 and the rest are upper neighbors.
Hence, $T$ is embedded as a rainbow subgraph.
\epf

It is important to note that this proposition does not imply
 $f(n,G|\cF)$ to be non-quadratic for every family $\cF$
 of trees.
Transparent counterexamples are the stars and the paths,
 as demonstrated by rainbow \tur\ graphs combined with
 monochromatic and \Lex\ colorings, respectively.
Tight asymptotics on these will be presented later.

\subsection{Large stars}

Next we consider two families $\cF$ defined in terms of stars
 parametrized with the fixed integer $r \geq 2$:
\begin{itemize}
 \item $\cF(r)$ --- the graphs without $K_{1,r}$ as a subgraph
  (i.e., $\Delta\leq r-1$);
 \item $\cF^*(r)$ --- the graphs without $K_{1,r}$ as an induced
  subgraph.
\end{itemize}
The family $\cF(r)$ is hereditary, and $\cF^*(r)$ is induced-hereditary.

The following observation will be useful.

\blm
If\/ $|G| \geq r\alpha(G) +1$, then every induced subgraph
 of\/ $G$ on\/ $r\alpha(G) +1$ vertices is non-\/$r$-colorable.    
\elm

\bpf
One of the basic facts in graph coloring is that the product
 $\chi\alpha$ is an upper bound on the order of the graph
 in question.
In the current case if $H\subset G$ is of order
 $r\alpha(G) +1$, then $\chi(G) \geq \chi(H)
 \geq \lceil\, |H|/\alpha(H) \rceil = r+1$.
\epf

Since all $(r-1)$-degenerate graphs are $r$-colorable, it
 follows that every non-$r$-colorable graph $G$ must contain an
 induced subgraph with minimum degree at least $r$.

\bpn
Suppose\/ $m \geq 2$ and\/
 $\lceil \, |G|/m \rceil \geq r\alpha(G) +1$.
Then, for any fixed\/ $r\geq 2$, both\/ $f(n,G|\cF(r))$ and\/
 $f(n,G|\cF^*(r))$ exceed\/ $\ex(n,K_{m+1})$.
\epn

\bpf
Consider the \rtdl$(n,m)$ pattern.
Some vertex class of its rainbow core $T_{n,m}$
 contains at least $r\alpha(G)+1$ vertices, hence by the
 preceding lemma a non-$r$-colorable $H \subset G$ occurs.
Consequently a graph $H^*\subset G$ with minimum degree
 at least $r$ is present in that class. 
Let $v$ be the top vertex of $H^*$ in the corresponding \Lex\ order.
Then $v$ has at least $r$ neighbors from $H^*$ in lower indices
 of \Lex\ and hence induces a monochromatic $K_{1,r'}$ edge
 class with $r'\geq r$ in $G$, which is not a member of $\cF(r)$.
\epf

This result shows that already with the partial information
 $|G|>r\alpha(G)$ a quadratic growth of $f(n,G|\cF)$ follows.
However, for an asymptotically tight formula one needs
 to know more details about $G$.
A result of this kind is Theorem 5.6 in \cite{ar-G}.

\subsection{Connectivity}

Requiring all members of $\cF$
 to be connected is a very strong condition with substantial
 impact on the behavior of $f(n,G|\cF)$.
We demonstrate this fact with two results.

\btm
Let\/ $\cF$ be obtained from the family of all connected graphs
 by deleting all complete graphs\/ $K_r$ with\/ $r\geq 3$.
Then, for every\/ $k\geq 4$ we have\/
 $f(n,K_k|\cF) = \ex(n,K_{k-1})+2$.
\etm

\bpf
To prove that $\ex(n,K_{k-1})+2$ is a lower bound
 on $f(n,K_k|\cF)$, we take the \rtem$(n,k-2)$ coloring.
Denote by $V_1,\dots,V_{k-2}$ the vertex classes of its rainbow
 $R\cong T_{n,k-2}$, and by $c$ the color assigned to the edges of $K_n-R$.
Then any copy of $K_k$ either has at least two vertices in each
 of two classes $V_i,V_j$, or has at least three vertices in a
 class $V_i$.
In the former case the color class $c$ induces a disconnected
 graph in $K_k$, which is not a member of $\cF$; and in the
 latter case it induces a $K_r$ with $r\geq 3$, which is
 not in $\cF$ either.

For the matching upper bound we observe that the inequality
 $f(n,K_k|\cF) \leq \AR(n,K_k)$ holds due to the definitions,
 and it is known that $\AR(n,K_k) = \ex(n,K_{k-1})+2$.
This was proved for large $n$ in \cite{ESS-75}, and for all
 $n\geq k$ in \cite{S-04}.
\epf

\btm
Let\/ $\cF$ be the family of all connected bipartite graphs.
If\/ $\chi(G)=k\geq 4$ and\/ $G$ contains an induced
 forest\/ $T$ such that\/ $\chi(G-T)=k-1$, then\/
 $\ex(n,K_{k-1})+2\leq f(n,G|\cF)\leq
  (1+o(1))\,\ex(n,K_{k-1})$.
\etm

\bpf
The lower bound is obtained in the same way as above:
 we take the \rtem$(n,k-2)$ coloring with rainbow
 $R\cong T_{n,k-2}$ and additional color $c$.
Then clearly $\chi(G) \leq \sum_{i=1}^{k-2} \chi(G_i)$.
Since $\chi(G)=k$, and $R$ has only $k-2$ vertex classes,
 any copy of $G$ in $K_n$ either contains
 edges from at least two classes of $R$, making the
 $c$-edge-class of $G$ disconnected, or contains a
 non-bipartite $c$-colored subgraph in one class of $R$.
Hence, in either case, color $c$ induces a subgraph of $G$
 which is not a member of $\cF$.

For the asymptotic upper bound we apply Theorem \ref{t:CCCL} and
 find a clean $(k-1)$-partite grapk $K\cong K_{t,\dots,t}$,
 where $t=|G|$ (in fact a very generous choice for $t$)
 and the edge coloring inside each vertex class is
 monochromatic or rainbow or \Lex.

Inside the first vertex class of $K$ we select a $T'\cong T$
 which is either monochromatic or rainbow, hence a member
 of $\cF$ in either case.
This is trivial to do if the class is monochromatic or rainbow;
 and can also be done via Proposition~\ref{p:tree-lex} if it is
 \Lex-colored.
And then $T'$ is extendable to an $\cF$-decomposed copy of $G$,
 using a proper number of vertices from each further class of
 $K_{t,\dots,t}$ to form the vertex classes of an $\cF$-colored
 copy of $G$.
\epf

\subsection{Disconnectedness}
\label{ss:disconn}

The notion considered here is complementary to the one in the
 previous subsection, but again quadratic growth will be observed.
Let us recall for $p\geq 5$ that $\lar(n,K_p) \sim
 \ex(n, K_{\lceil p/2 \rceil})$, which means the family
 $\cM:=\{mK_2 : m\geq 1\}$ of matchings in the terminology
 of the $f$ function.
It turns out that quadratic growth order remains the same if
 $\cM$ is supplemented with all disconnected graphs, but on the
 other hand the same asymptotics remain valid if just $K_2$ is
 taken together with a single suitably chosen disconnected graph.
Note that the family $\cfs$ below contains $\cM$.

\btm
Let\/ $p,s$ be integers where\/ $p\geq 5$ and\/ $2\leq s\leq p-2$.
Further, let\/ $\cF$ be the family of all disconnected graphs
 without isolated vertices, and set\/ $\cfs = \cF \cup \{K_2\}$.
Then both\/ $f(n,K_p|\cfs)$ and\/ $f(n,K_p|\{K_2,K_s\cup K_{p-s}\})$
 are at least\/ $\ex(n, K_{\ceil{p/2}}) + \ceil{p/2}$
  and at most\/ $(1+o(1))\,\ex(n, K_{\lceil p/2 \rceil})$.
\etm

\bpf
Since both $K_2$ and $K_s \cup K_{p-s}$ are members of $\cfs$, by
 Proposition \ref{p:monoton}~$(i)$ it suffices to prove the lower
 bound for $\cfs$ and the upper bound for $\{K_2,K_s\cup K_{p-s}\}$.

For $\cfs$ we can take the \rtdm$(n,\lceil p/2 \rceil -1)$ coloring.
Then any copy of $K_p$ contains at least three vertices in
 some vertex class of the rainbow $T_{n,\lceil p/2 \rceil -1}$,
 inducing a color which defines a connected (complete) subgraph
 in $K_p$ and therefore not a member of $\cfs$.
Hence, more than $\ex(n, K_{\ceil{p/2}}) + \ceil{p/2}-1$ colors
 are needed to force the presence of a suitable copy of $K_p$.

For $\{K_2,K_s\cup K_{p-s}\}$ we apply Theorem \ref{t:CCCL}.
If more than $\ex(n, K_{\lceil p/2 \rceil}) + h(n)$ colors are
 used, we find a clean rainbow $K_{p,\dots,p}$ with $r=\lceil p/2 \rceil$
 homogeneously colored vertex classes $V_1,\dots,V_r$ of size $p$ each.
If at least one class $V_i$ is rainbow, we take it as $K_p$.
If at least two classes $V_i,V_j$ are monochromatic in the same
 color, we take $s$ vertices from $V_i$ and $p-s$ vertices from
  $V_j$; those $p$ vertices induce a required copy of $K_p$,
 in which one color class is $K_s\cup K_{p-s}$ and the
 other edges have their private colors.
In any other case we can take a rainbow matching of size
 $\lceil p/2 \rceil$ whose edges are inside mutually distinct
 vertex classes.
Since $K_{p,\dots,p}$ is rainbow and clean, this yields a rainbow
 copy of $K_p$.
\epf

The other side of the coin is when only matchings of relatively
 small sizes are allowed, nothing else.
The result for this sparse small family $\cfs_r$ gives a smooth
 transition from the local \arr\ number $\lar(n,K_p)\sim
 \ex(n,K_{\lceil p/2 \rceil})$ to the (standard) \arr\ number
 $\AR(n,K_p)\sim \ex(n,K_{p-1})$.

\btm   \label{t:r-match}
Let\/ $r\geq 1$ be an integer, and let\/ $\cfs_r = \{mK_2 \mid
 1\leq m\leq r\}$ be the family of matchings of size at most\/ $r$.
If\/ $r\geq 2$ and\/ $p>2r$, or\/ $r=1$ and\/ $p\geq 4$, then\/
 $$\ex(n, K_{p-r}) + 2
 \leq f(n, K_ p |\cfs_r) \leq (1+o(1))\,\ex(n, K_{p-r}) \,.$$
The same bounds are valid if\/ $\cfs_r$ is supplemented with all
 non-matching triangle-free graphs (connected and disconnected
 as well).
\etm

\bpf
For the lower bound we now take the \rtem$(n,p-r-1)$ color pattern.
In this case a complete subgraph satisfying the conditions of $g$
 cannot meet a vertex class of the rainbow core $T_{n,p-r-1}$ in
 more than two vertices, because a monochromatic $K_s$ with $s>2$
 is not allowed in $\cfs_r$, nor in any triangle-free graph.
Hence, the largest complete subgraph having its color partition
 with members of $\cfs_r$ and/or non-matching
  triangle-free graphs can have at
 most $p-1$ vertices, because it can meet only $r$ classes of the
 rainbow $T_{n,p-r-1}$ in two vertices, and contain only one
 vertex from each of the other $p-2r-1$ classes.
Hence, we need at least $\ex(n, K_{p-r}) + 2$ colors in order to
 guarantee a $K_p$ subgraph compatible with $\cfs_r$ for sure.

For the upper bound we apply Theorem \ref{t:CCCL} for the
 particular case of $K\cong K_{2,\dots,2}$ with $p-r$ vertex classes.
As it has been proved, to obtain a clean rainbow $K$ it suffices to have
 no more than $\ex(n,K_{p-r})+h(n)$ colors, where $h(n)=o(n^2)$.
Keeping $r$ vertex classes of $K$, and selecting one vertex from
 each of its other $p-2r$ vertex classes, an $\cfs_r$-colored
 copy of $K_p$ is obtained, regardless of the color distribution on
 the edges contained in the $r$ selected vertex classes of~$K$.
\epf

We conclude this section with a finite mixture of connected and
 disconnected graphs, mentioned already at the end of
 Section \ref{ss:stab}.

\btm   \label{t:k-small}
Let\/ $k \geq 2$ be any integer and\/ $\cF(k)$ the family of
 all graphs on at most\/ $k$ vertices.
Then for every\/ $p \geq k+2$ we have\/ $f(n,K_p|\cF(k)) =
 ( 1+o(1))\,\ex(n,K_{p-k+1})$.
\etm

\bpf
For a lower bound we take the \rtem$(n,p-k)$ coloring of $K_n$.
In a complete subgraph the edges from the big color class can
 span at most $k$ vertices.
Even when they belong to the largest vertex class of the rainbow
 $T_{p-k}\subset K_n$, only one additional vertex can be
 taken from each of the other $p-k-1$ vertex classes.
Hence, no $K_p$ compatible with $\cF(k)$ exists.

For an upper bound we apply Theorem \ref{t:CCCL}.
In any edge coloring of $K_n$ with more than $\ex(n,K_{p-k+1})+h(n)$ colors
 we find a clean rainbow $K=K_{k,1,\dots,1}\cong K_p-K_k$.
No matter how the big vertex class is edge-colored, it is
 compatible with $\cF(k)$, because its colors do not appear in $K$.
\epf

\section{Stability concerning chromatic number}
\label{ss:stab}

In \cite{ar-G}, a graph was defined to be stable with respect
 to a hereditary family
  $\cF$, if $\chiF(G) = \chi(G)$; that is, if $\chi(G - D)
  = \chi(G)$ holds whenever $D\subset G$ and $D\in \cF$.
This notion allows us to prove that any (not necessarily
 hereditary) $\cF$ allows infinitely
 many graphs $G$ to have $f(n,G|\cF)$ with quadratic growth,
 unless $\cF$ contains graphs that are very dense.

\btm   \label{t:k2n2}
Let\/ $k\geq 3$ be any fixed integer, and\/ $\cF$ any family
 of graphs that contains\/ $K_2$.
If there is a threshold\/ $n_0$ such that\/
 $\max\{ e(D) : D\in \cF, \, |V(D)|=n\} < n^2/k^2$ holds
 for every\/ $n>n_0$, then there exist infinitely many graphs\/
 $G$ with\/ $\chiF(G)=\chi(G)=k$, and all those\/ $G$ have\/
 $\ex(n,K_{k})+2\leq f(n,G|\cF)\leq
  (1+o(1))\,\ex(n,K_{k})$.
\etm

\bpf
We prove that the complete $k$-partite graphs $K=K_{p,\dots,p}$
 satisfy the conditions for all $p$ sufficiently large.
Consider any $D\in \cF$.
Each edge in any copy of $D$ in $K$ is contained in $p^{k-2}$
 cliques $K_k$ of $K$.
By assumption, if $kp>n_0$ then any $D\in \cF$ of order at most
 $kp$ can have fewer than $(kp)^2\!/k^2=p^2$ edges. 
Hence, removing all the edges of any $D$, only fewer than $p^k$
 copies of $K_k$ can be destroyed.
Consequently $K_k\subset K-D$ and $\chiF(K)=k$ are valid.

Concerning the inequalities on $f(n,G|\cF)$, for the lower bound
 we can take the \rtem$(n,k-1)$ pattern.
Recall that here the rainbow $T_{n,k-1}$ has a
 monochromatic complement.
Since the large color class of $E(K_n)$ can accommodate only one
 color class of $G$, say $F\in\cF$, the remaining $k$-chromatic
 $G-F$ should be part of the $(k-1)$-chromatic rainbow $T_{n,k-1}$, which is impossible.
On the other hand, since $K_2\in\cF$ and $\chi(G)=k$, we have
  $$
    f(n,G|\cF) \leq \AR(n,G) \leq \ex(n,K_{|G|,|G|,\dots,|G|}) + 1
     \leq (1+o(1))\,\ex(n,K_{k})$$
  by the \erd--Stone--Simonovits
 theorem, hence proving the claimed upper bound.
\epf

An interesting example is the family $\cF$ of stars.
It is hereditary, hence the methods developed in \cite{ar-G}
 perfectly handle it.
Moreover it is sparse, with $e(D) = |V(D)|-1$.
On the other hand it is worth observing that the condition
 $\chiF(G)=\chi(G)=k$ holds if and only if $G$ is a
 $k$-chromatic graph which is not vertex-critical.
Further, with respect to rainbow coloring the condition
 ``vertex-critical'' is replaced with ``edge-critical''.

\bex
A simple natural case is provided by the collection of all graphs
 on $k$ vertices (\/$k\geq 3$),
 to which we shall return in Section \ref{ss:disconn}.
Also with respect to this family, there exist infinitely many
 graphs satisfying the stability property described in
 Theorem \ref{t:k2n2}.
\eex

The following very rich class of families may be viewed as a
 generic case for which the principle of stability applies. 

\bex
Let\/ $H$ be any bipartite graph.
Consider the vamily\/ $\cF$ of all graphs not containing\/ $H$
 as a subgraph.
Since\/ $H$ is bipartite, we have\/ $\ex(n,H)=o(n^2)$, hence the
 conditions of Theorem \ref{t:k2n2} hold true for every\/ $k$ and
 every\/ $n>n_0(H,k)$.
\eex

\section{Excluding just one graph: Complement of exact piercing}
\label{s:excl1}

In this section we consider the most flexible cases possible,
 namely where the color classes can induce anything but just one
 excluded graph.
These results, when combined with the complementary property
 (Proposition \ref{p:compl}~$(ii)$), determine exact piercing
 $g(n,K_p|G)$ tightly if $G$ is not a star, and asymptotically
 if $G$ is a star $K_{1,t}$, for $t \geq 2$ and $p \geq 2t +1$.

Formally, for a graph $H\neq K_2$, denote by $\nem{H}$ the
 family of all finite graphs without isolates, except $H$.
It turns out that there is a substantial distinction whether
 $H$ is a star or a non-star graph.

\btm   \label{t:excl-star}
\nev{Excluding one star} \
Let\/ $p > t \geq  2$ be any integers.
\begin{itemize}
 \item[$(i)$] If\/ $p \geq  2t+1$, then\/
 $f(n,K_p|\nem{K_{1,t})) = (1+o(1))\,\ex(n,K_{\ceil{p/t}})$.
 \item[$(ii)$] If\/ $t+1 \leq p \leq  2t$, then\/
  $f(n,K_p |\nem{K_{1,t}}) \leq O(n^{2-1/t})$.
  Moreover, for all\/ $n\geq p\geq t+1$,
   $f(n,K_p |\nem{K_{1,t}}) \geq (p-t-1)(n-p+t+1) + \binom{p-t-1}2} + 1$.
\end{itemize} 
\etm

\bpf
$(i)$ \
The lower bound $f(n,K_p|\nem{K_{1,t}}) \geq \ex(n,K_{\ceil{p/t}})
 + n - \ceil{p/t} + 1$
 is obtained from the composite color pattern \rtdl$(n,\ceil{p/t}-1)$.
Indeed, among any $p$ vertices at least $t+1$ belong to the
 same vertex class of $T_{n,\ceil{p/t}-1}$, and then \Lex\
 induces a forbidden monochromatic class $K_{1,t}$.

The upper bound
 $f(n,K_p|\nem{K_{1,t}}) = (1+o(1))\,\ex(n,K_{\ceil{p/t}})$,
 that matches the lower bound in its asymptotics,
 follows by Theorem \ref{t:CCCL}.
In any coloring with that many colors, we can find a clean
 complete $\ceil{p/t}$-partite rainbow subgraph in which all
 vertex classes have size $t$, hence leaving no room for a
 monochromatic $K_{1,t}$.

\ssk

\nin
$(ii)$ \
By the same method, applying Theorem \ref{t:CCCL}, we find a
 clean rainbow $K_{t,t}$ subgraph; to do so, we do not need
 more than $O(n^{2-1/t})$ colors.
Then take a whole vertex class of $t$ vertices, and supplement it
 with $p-t\leq t$ vertices from the other class.
In the obtained $K_p$ subgraph none of the color classes can be
 a star $K_{1,t}$.

The claimed lower bound can be verified by taking a monochromatic
 $K_{p-t-1}$, and assigning a private color to each edge of
 $K_n-K_{p-t-1}$.
Then from any $p$ vertices, at least $t+1$ belong to $K_{p-t-1}$,
 and that color class is a complete graph, not $K_{1,t}$.
Every other color class is a single edge.
\epf

One may observe that in the range $t+1\leq p \leq 2t-2$
 it suffices to find a rainbow $K_{p-t,t}$ and the strength
 of Theorem \ref{t:CCCL} is not needed;
 the result then follows just by the fact
 $\ex(n,K_{t-2,t})=o(n^2)$.

\btm   \label{t:excl-nonstar}
\nev{Excluding one non-star graph} \
If\/ $G$ is not a star, and
  \begin{itemize}
   \item[$(i)$] either\/ $G$ is non-complete and\/ $p\geq |G|$,
   \item[$(ii)$] or\/ $G$ is complete and\/ $p > |G|$,
  \end{itemize}
 then\/ $f (n,K_p |\nem{G}) = 1$ holds for every\/
 $n$ sufficiently large with respect to\/~$p$.
\etm

\bpf
Consider any edge coloring of $K_n$.
According to Theorem \ref{t:E--R}, for every $n\geq n_0(p)$,
 there occurs a $K_p$ which is either monochromatic or rainbow
 or \Lex-colored.
Under the current assumptions, none of these patterns can induce
 a color class isomorphic to $G$.
(If $G$ is complete, then $p>|G|$, hence if the copy of $K_p$
 is monochromatic, then $G$ is only a proper subgraph of
 that color class.)
\epf

\section{Quadratic growth via integer sequences}
\label{s:quad-spec}

Here we consider some rather peculiar classes, for which our
 methods still allow to prove tight asymptotic results,
 despite that those classes are very far from being hereditary.
The common principle in forming those classes is that we take
 a finite or infinite sequence $S$ and consider the family
 $\fs{S}$ of all graphs whose number of edges is in $S$. 

In the lower-bound constructions of this section it will suffice
 to consider the \rtdl$(n,m)$ colorings with suitably chosen $m$.

\subsection{Graphs with square numbers of edges}

Let us denote by $\fsq$ the family of all graphs having a square
 number of edges,
 $\fsq = \{ G \mid e(G) = k^2, \ k\geq 1 \}$.

\btm   \label{t:square}
For every\/ $p\geq 5$,
 $f(n,K_p|\fsq) = (1+o(1))\,\ex(n,K_{\ceil{p/2}})$.
\etm

\bpf
For the lower bound
 $f(n,K_p|\fsq) \geq \ex(n,K_{\ceil{p/2}}) + n - \ceil{p/2} + 2$,
 we take the \rtdl$(n,\ceil{p/2} - 1)$ coloring.
The number of colors is $\ex(n,K_{\ceil{p/2}}) + n - \ceil{p/2} + 1$.
From any $p$ vertices, at least three belong to the same
 vertex class of the rainbow $T_{n,\ceil{p/2} - 1}$,
 and then a color class $P_3$ with just two edges
 occurs, which is not allowed in $\fsq$.
(The color of that $P_3$ does not occur in any other vertex class.)

We derive the upper bound $f(n,K_p|\fsq) \leq
 (1+o(1))\,\ex(n,K_{\ceil{p/2}})$ on applying Theorem \ref{t:CCCL}.
With the claimed number of colors, we can guarantee a complete
 $\ceil{p/2}$-partite rainbow subgraph $K$ with $p$ vertices
 in every partite set, and each of those sets induces a
 monochromatic or rainbow or \Lex-colored $K_p$.

If a rainbow vertex class arises, then in the $K_p$ induced
 by it, every color class is a $K_2\in\fsq$.
If at most one part is monochromatic, or if the monochromatic
 parts have mutually distinct colors, then we select one edge
 from each partite set of $K$, first from the monochromatic ones, and
 then from the \Lex\ parts, in a way that the selected matching
 is rainbow.
This is possible because each \Lex\ part offers $p-1$ colors
 for choice.

Finally, if there are two monochromatic parts $V_i,V_j$ of the
 same color, let us write $p$ in the form $p=2q-1+\epsilon$,
 where $q=\ceil{p/2}$ (hence $\epsilon=0$ if $p$ is odd, and
  $\epsilon=1$ if $p$ is even).
We compose $K_p\subset K_n$ with $q$ vertices from $V_i$,
 $q-1$ vertices from $V_j$, and $\epsilon$ vertices from a
 third vertex class.
Then, in this $K_p$, each edge connecting distinct parts of $K$
 forms a color class $K_2\in\fsq$, and the other edges form a
 monochromatic $K_q\cup K_{q-1}\in\fsq$ with
 $\binom{q}{2} + \binom{q-1}{2} = (q-1)^2$ edges.
Note that if $\epsilon=1$, then the last vertex can be selected
 because the number of parts in $K$ is $\ceil{p/2}\geq 3$,
 as $p\geq 5$.
\epf

\subsection{Arithmetic progressions}

Let $k\geq 2$ be any integer, and consider the family $\fek$
 of all graphs $G$ with $e(G) \equiv 1$ (mod $k$).
In particular, $\cF_{1,2}$ consists of all graphs whose number
 of edges is odd.

\btm   \label{t:1-mod-k}
For every\/ $p\geq 2k-1$ if\/ $k$ is odd, and\/ $p\geq 4k-1$ if\/ $k$ is even,\/
 $f(n,K_p|\fek) = (1+o(1))\,\ex(n,K_{\ceil{p/2}})$.
\etm

\bpf
Also here, the lower bound is obtained from the
 \rtdl$(n,\ceil{p/2}-1)$ coloring.
Then every $K_p\subset K_n$ induces a monochromatic
 $P_3\notin\fek$ as a color class.
As a matter of fact, this construction works for all $p\geq 5$,
 independently of the value of $k$.

For the upper bound we choose an integer $t$ such that
 $p \geq t \geq \ceil{(p+1)/2}$ and $\binom{t}{2} \equiv 1$ (mod $k$).
Under the assumptions put on $p$, such a $t$ always exists:
 \begin{itemize}
  \item If $k$ is odd, let $t\equiv k-1$ (mod $k$).
   Then $t(t-1)/2 \equiv (k-1)(k-2)/2 = (k^2 - 3k + 2)/2
 = k(k-3)/2 +1 \equiv 1$ (mod $k$),
  because $(k-3)/2$ is an integer for $k$ odd. 
  \item If $k$ is even, let $t \equiv 2k-1$ (mod $2k$).
   Then $t(t-1)/2 \equiv (2k-1)(2k-2)/2 =
    (4k^2 - 6k +2)/2 = 2k^2 - 3k +1 \equiv 1$ (mod $k$).
 \end{itemize}
Observe that the interval $[\,\ceil{(p+1)/2},p]$ contains exactly
 $\floor{(p+1)/2}$ consecutive integers, which means, under the given assumptions
 on $p$, at least $k$ for $k$ odd and at least $2k$ for $k$ even.
Thus, a required $t$ exists in every case.
Moreover, it is important to note that $(\ceil{p/2}-1)+t \geq p$
 holds because $t \geq \ceil{(p+1)/2}$.

Now, if the number of colors in $K_n$ is sufficiently large,
 we can apply Theorem \ref{t:CCCL} to obtain a complete
 $\ceil{p/2}$-partite, clean, rainbow subgraph $K\subset K_n$,
 with $t$ vertices in each partite set, each part being
 monochromatic or rainbow or \Lex-colored.

If a rainbow part occurs, then we take its $t$ vertices and
 supplement them with one further vertex from $p-t \leq
 \ceil{p/2}-1$ other parts, hence obtaining a rainbow $K_p$.
If a monochromatic part occurs, we do the same, obtaining a
 $K_p$ in which a monochromatic $K_t$ has $\binom{t}{2}\equiv 1$
 (mod $k$) edges, and all edges in this $K_p-K_t$ have their
 private color.
Finally, if each part is \Lex-colored, we take a rainbow
 matching by selecting one edge from each part.
In this way a rainbow $K_{2\ceil{p/2}} \supseteq K_p$ is obtained.
\epf

In case of $k=2$, that is for the family of graphs with odd
 numbers of edges, it can be proved that the general asymptotic
 formula is valid in the entire range $p\geq 5$.

\btm
For every\/ $p\geq 5$ we have\/
 $f(n,K_p|\cF_{1,2}) = (1+o(1))\,\ex(n,K_{\ceil{p/2}})$.
\etm

\bpf
Our general Theorem \ref{t:1-mod-k} settles all cases $p\geq 7$,
 so it remains to consider $p=5$ and $p=6$ only.
Then $\ceil{p/2}=3$, and the \rtdl$(n,2)$ coloring provides us
 with a suitable lower bound, as it does not admit
 any $K_p$ subgraph without a 2-edge $P_3$ color class.

Turning to the asymptotic upper bound, for $p=5$ we can again
 refer to the preceding proof.
Namely, in this case the value $t=3=\ceil{5+1}/2$ is a suitable
 choice, as it luckily is of the form $2k-1$.
But this method does not work for $p=6$, because the integers
 of the interval $[4,6]$ do not satisfy the needed properties.

Instead, for $p=6$, we start with a clean rainbow $K_{5,5,5}$
 (rather than $K_{3,3,3}$) with homogeneous patterns in its
  parts, found with the help of Theorem \ref{t:CCCL}.
If a rainbow part occurs, we supplement it with one vertex
 from another part and obtain a rainbow $K_6$.
Also, if at least two parts are \Lex-colored, we can find a
 rainbow $K_6$ by selecting a rainbow matching $3K_2$ from the
 edges of the three parts (starting with an edge from the
 monochromatic part, if not all parts are \Lex).
We can proceed in a similar way if there are two or three
 monochromatic parts of mutually distinct colors.

The selection of one edge from each part works also if we have
 three monochromatic parts of the same color.
Then the found $K_6$ has one color class $3K_2$, and all its
 other edges have their private color.
Hence, we are left with the case where two monochromatic parts,
 say $V_1$ and $V_2$, have the same color, and the third part
 $V_3$ is either \Lex-colored or monochromatic in another color.
Then it is a proper choice to select a $K_3$ from $V_1$,
 one vertex from $V_2$, and an edge from $V_3$, whose color is
 distinct from the color of $V_1$.
In this way $K_3$ is a 3-edge color class, and
 every edge of $K_6-K_3$ has its private color.
\epf

\subsection{Integer sequences with less regularity}

After some explicitly defined sequenes (the squares, and the
 numbers whose residue modulo $k$ equals 1), here we consider 
 the following more flexible property.

\bdf   \nev{Legal sequence, S-good integer}
We say that a finite or infinite sequence\/ $S$ of positive
 integers is legal if\/ $1\in S$ and\/ $2\notin S$.
Further, given a legal sequence\/ $S$, an integer\/ $p \geq 5$ is
 called\/ $S$-good if there is an integer\/ $t$ such that\/
 $p \geq t \geq \ceil{(p+1)/2}$ and\/ $\binom{t}{2} \in S$.
\edf

So, in connection with $S$-good integers, one may restrict
 attention to sequences of triangular numbers.

\bex   \nev{Special sets of triangular numbers}
 \begin{enumerate}
  \item According to \cite{M-89}, the set of triangular
    Fibonacci numbers is\/ $S=\{1,3,21,55\}$.
   The corresponding\/ $t$-values are\/ $t=2,3,7,11$.
   Then:
    \begin{itemize}
     \item For\/ $t= 2$, no\/ $p \geq 5$ is\/ $S$-good.
     \item For\/ $t= 3$, $p = 5$ is\/ $S$-good.
     \item For\/ $t= 7$, $p$ is\/ $S$-good if\/ $7 \leq p \leq 13$.
     \item For\/ $t= 11$, $p$ is $S$-good if\/ $11 \leq p \leq 21$.
    \end{itemize}
  \item According to \cite{M-91}, the set of triangular
    Lucas numbers is\/ $S=\{1,3,5778\}$.
   The corresponding\/ $t$-values are\/ $t=2,3,108$.
    \begin{itemize}
     \item For\/ $t= 108$, $p$ is\/ $S$-good if\/ $108 \leq p \leq 215$.
    \end{itemize}
  \item According to \cite{Argy}, the set of triangular
    Mersenne numbers\/ $2^n-1$ is\/ $S=\{1,3,15,4095\}$.
   The corresponding\/ $t$-values are\/ $t=2,3,6,91$.
    \begin{itemize}
     \item For $t= 6$, $p$ is $S$-good if $6 \leq p \leq 11$.
     \item For $t= 91$, $p$ is $S$-good if $91 \leq p \leq 181$.
    \end{itemize}
 \end{enumerate}
\eex

The relevance of legal sequences and $S$-good integers is
 demonstrated by the following result.

\btm   \label{t:S-good}
Let\/ $S$ be a legal sequence and\/ $p\geq 5$ an\/ $S$-good integer.
Then\/ $f(n, K_p |\fs{S}) = (1+o(1))\,\ex(n,K_{\ceil{p/2}})$.
\etm

\bpf
As we have already seen, the \rtdl$(n,\ceil{p/2}-1)$
 coloring satisfies the property that
 every $K_p$ subgraph induces a $P_3$ color class.
This class has just two edges, hence not allowed for $f$ because
 $S$ is legal, $2\notin S$.

For the upper bound we apply Theorem \ref{t:CCCL} to find a clean,
 rainbow, complete $\ceil{p/2}$-partite graph with homogeneous parts
 of size $p$.
If a rainbow part occurs, we take it as $K_p$.
If all parts are \Lex-colored, we choose a rainbow matching of
 $\ceil{p/2}$ edges.
If there is a monochromatic part, we take a $K_t$ inside this
 part, and supplement it with one vertex from $p-t$ other parts.
This last step is also feasible because $p-t\leq \ceil{p/2}-1$.
It yields one color class where the number of edges is
 $\binom{t}{2}\in\fs{S}$, and all the other color classes are
 single edges.
\epf

\bex
  Let\/ $S=\{\mathrm{triangular~numbers}~\binom{q}{2}~\mathrm{of~primes}~q\}$.
   Then\/ $1 \in S$ by\/ $q = 2$, and\/ $2$ is non-triangular,
     hence\/ $2 \notin S$.
    Furthermore, inverting Bertrand's postulate, for every\/
    $p \geq 5$ there is a prime\/ $t$ such that\/
     $p \geq t \geq \ceil{(p+1)/2}$;
     e.g.,\/ $t=3$ for\/ $p = 5$,\/ $t = 5$ for\/ $p = 6$, etc.
   Thus:
 \begin{itemize}
     \item Every integer\/ $p \geq 5$ is\/ $S$-good.
 \end{itemize}
\eex

This example motivates the following notion, which is
 supplementary to the ``$S$-good'' property.

\bdf   \nev{Well-spaced sequence}
An infinite sequence\/ $S$ of triangular numbers is well-spaced if\/
 $1 \in S$,\/ $2 \notin S$, and for every\/ $p \geq 5$ there is
 an integer\/ $t$ such that\/ $p \geq t \geq \ceil{(p+1)/2}$ and\/
  $\binom{t}{2} \in S$.
\edf

If $t\geq 3$ is an integer with $\binom{t}{2}\in S$, then the values
 $p\geq 5$ which are $S$-good due to $t$ constitute the interval
  $[t',2t-1]$,  where $t'=t$ if $t\geq 5$ and $t'=5$ if $t=3,4$.
Based on this reason, the following assertion can be posed.

\bpn
The sparsest well-spaced sequence\/ $S=\{s_i \mid i\geq 1\}$ has
 the elements\/ $1,10,45,190,\dots$ generated by the sequence\/
 $t_1=2$ and\/ $t_i=5\cdot 2^{i-2}$ for\/ $i\geq 2$
  (or,\/ $\floor{5\cdot 2^{i-2}}$ for all\/ $i\geq 1$); hence
  $$
    s_1 = 1 \qquad \mathrm{and} \qquad
     s_i = \binom{5\cdot 2^{i-2}}{2} \quad \mathrm{for~all}
     \quad i\geq 2 \,,
  $$
 with its\/ $t$-sequence\/ $t_1,t_2,t_3,t_4,\ldots = 2, 5, 10, 20,\dots$.
\epn

\bpf
If $S$ is well-spaced and $s_i=\binom{t_i}{2}$, where
 the corresponding $t$-sequence is
 $t_1<t_2<\cdots$\,, then we must have $t_1=2$ to guarantee
  $1 \in S$, and $3\leq t_2\leq 5$ to make $p=5$ $S$-good.
Since any $t_i\geq 3$ settles the $S$-good status of integers
 at most $2t_i-1$, to have $S$ well-spaced and ensure that
  $p=2t_i$ is $S$-good, we need $t_{i+1}\leq 2t_i$.
So, in general, we have $t_{j+2}\leq 5\cdot 2^j$ for all $j\geq 0$.
On the other hand, it is readily seen that the claimed sequence
 is well-spaced.
\epf

\btm   \label{t:well-sp}
If\/ $S$ is a well-spaced infinite sequence of integers,
 then\/ $f(n, K_p |\fs{S}) = (1+o(1))\,\ex(n,K_{\ceil{p/2}})$
 holds for all\/ $p\geq 5$.
\etm

\bpf
If $S$ is well-spaced, then every $p\geq 5$ is $S$-good, and
 Theorem \ref{t:S-good} applies.
\epf

\section{Towards subquadratic growth}   \label{s:subqu}

Certainly, if $\cF$ is the family of all graphs, then we have
 $f(n,G|\cF) = 1$ for all $G$ and all $n \geq |G|$.
Further, if $K_2$ is not a member of $F$, then $f(n,G|\cF)$
 doesn’t exists, as shown by the rainbow coloring of $K_n$. 
On the other hand, we have seen many classes of families $\cF$
 for which $f(n,G|\cF)$ is a quadratic function of $n$.
In this context the results of this short section give
 contributions to the following natural problem.

\bpm
Characterize the families\/ $\cF$ minimal under inclusion,
 for which\/ $f(n,G|\cF)$ is subquadratic for every\/ $G$.
\epm

Based on observations from erlier sections, we can
 propose the following approach.

\bpn   \label{p:no-forest}
Let\/ $\cF$ be obtained from the family of all graphs by
 omitting all forests except stars.
Equivalently,\/ $\cF$ is the union of
 \begin{itemize}
  \item the family of all graphs containing a cycle, and
  \item the family of all stars.
 \end{itemize}
Then the following holds.
 \begin{itemize}
  \item[$(i)$] If\/ $G$ is acyclic, then\/ $f(n,G|\cF)=O(n)$.
  \item[$(ii)$] Otherwise\/ $f(n,G|\cF)=1$ for all\/ $n$
   large enough.
 \end{itemize}
\epn

\bpf
If $G$ is an acyclic graph, then, as $K_2\in\cF$, we can
 apply the general inequality chain
 $$
   f(n,G|\cF) \leq \AR(n,G) \leq \ex(n,G) +1 = O(n)
 $$
  valid for all graphs except the very last step, which is
 now applicable under the assumption that $G$ is a forest,
 hence its \tur\ number is linear in $n$.
Otherwise, if $G$ contains a cycle, we assume that $n$ is large
 and apply Theorem \ref{t:E--R}. 
Then any edge coloring of $K_n$ contains a $K_{|G|}$ which is
 either monochromatic or rainbow or \Lex-colored.
Monochromatic $G$ is in $\cF$ since $G$ has a cycle.
In rainbow $G$ each color class is $K_2$, a member of $\cF$.
In \Lex-colored $G$ all color classes are stars, which are
 members of $\cF$.
\epf

A further result offers a possibility to verify growth rates
 in the entire range between constant and quadratic:

\bpn   \label{p:small-AR}
Let\/ $h(n)$ be a positive integer function.
Let\/ $\cF$ be obtained from the family of all graphs by
 omitting all\/ $H$ such that\/ $\AR(n,H) \leq h(n)$
 and\/ $H$ is not a star.
 \begin{itemize}
  \item[$(i)$] If\/ $G \notin \cF$, then\/
   $f(n,G |\cF)\leq h(n)$.
  \item[$(ii)$] If\/ $G \in \cF$, then\/ $f(n,G|\cF)=1$
   for all large\/ $n$.
 \end{itemize}
\epn

\bpf
If $G \notin \cF$ then $f(n,G |\cF)\leq \AR(n,G) \leq h(n)$
 holds by definition.     
So, assume $G \in \cF$. 
According to Theorem \ref{t:E--R}, for $n$ large enough,
 with any number of colors in $K_n$ there is either a
 monochromatic copy of $G$ which is a member of $\cF$,
 or rainbow copy of $G$ whose color classes are single edges
 hence members of $\cF$, or \Lex-colored whose classes
 decompose $G$ into stars which are members of $\cF$.
\epf

\section{The piercing number}
\label{s:pierc}

In this section we study the function $\gnf$.
Recall that $K_2\notin\cF$ is assumed in the definition of $\gnf$.
As a preliminary observation, we state the following simple inequality.

\bpn   \label{p:UB-AR}
For every graph\/ $G$ and family\/ $\cF$, we have\/ $\gnf < \AR(n,G)$.
\epn

\bpf
In every coloring with at least $\AR(n,G)$ there occurs a
 rainbow copy of $G$.
Since $K_2\notin\cF$, this copy violates the conditins on $g$.
\epf

\subsection{Non-quadratic piercing}

The next result presents a class of instances which
 does not admit quadratic growth.

\btm   \label{t:connbip}
Let\/ $H$ be a connected bipartite graph of order\/
 $|H| = t+1 \geq 4$.
If\/ $2t \geq p \geq t +1$, then\/
 $g(n, K_p|\upp{H}) = O(n^{2 - \epsilon})$ for some\/
 $\epsilon = \epsilon(H) > 0$.
In particular,\/ $g(n, K_p|\upp{P_q}) = o(n^2)$ and\/
 $g(n, K_p|\upp{K_{1,q-1}}) = o(n^2)$ whenever\/
 $p\leq 2q-2$.
\etm

\bpf
Assume first that $H$ is 2-edge-connected.
In this case we apply the fact that
 $\AR(n, K_{t ,t}) \leq \ex(n,K_{t,t}) \leq O(n^{2-1/t})$,
 hence $\epsilon(H) = 1/t$ is a proper choice.
If an edge coloring uses more than that many colors, a
 rainbow $K'\cong K_{t,t}$ occurs, say with vertex classes $A,B$.
Let then $K\cong K_p$, $V(K)\subset V(K')$, such that the
 edge set of $H$ can be packed inside $E(K)\cap E(K')$.
(That is, the $p$ vertices of $K$ are suitably distributed
 between $A$ and $B$.)
This can be done, since $|H|=t+1\leq p\leq 2t$.

There is not enough space for a copy of $H$ inside $A$ or $B$,
 hence a monochromatic copy of $H$ should contain vertices
 in both parts of the rainbow $K'$.
But 2-edge-connectivity requires at least two edges between
 $A$ and $B$, thus no monochromatic $H$ can occur.

Assume now that $H$ is connected but not 2-edge-connected, i.e.,
 it contains a bridge.
Using Theorem \ref{t:CCCL} we find a clean copy $K'$ of
 $K_{t,t}$; and choose $K\cong K_p$ as above.
If $H'\cong H$ and $H'\subset K$, both sides of $K'$ must contain
 vertices of $H'$.
Now, since $H$ is connected, there is an edge
 $e\in E(K')\cap E(H')$, and the color of $e$ is not repeated
 inside $K$ because $K'$ is rainbow and clean.
Consequently, $K$ does not contain any monochromatic subgraph
 from $\upp{H}$, and the upper bound follows.
\epf

We note that the assumptions of the theorem cannot be improved.
If $p> 2t$ is allowed, then even dropping the requirement that
 $H$ is bipartite, we shall see that already the densest case
  $H=K_{t+1}$ yields quadratic growth for $g(n,K_p|\upp{H})$
 (and then of course every subgraph of $K_{t+1}$ does the same).
Also, connectivity is substantial: If $t=2s-1$ and $H=sK_2$,
 let $p=t+2=2s+1$.
Then $g(n,K_p|\upp{H})>\ex(n,K_3)=\floor{n^2\!/4}$.

\subsection{Piercing by small graphs}

In the following, we determine the exact values for some further
 families $\upp{F}$ generated by small graphs $F$.

The following example shows that the upper bound
 in Proposition \ref{p:UB-AR} is tight.

\bpn
For every\/ $n$ we have\/ $g(n,K_3|P_3)=g(n,K_3|\upp{P_3})=n-1$.
\epn

\bpf
Since $g(n,K_3|P_3)\leq g(n,K_3|\upp{P_3})$, it suffices to prove
 $g(n,K_3|P_3)\geq n-1$ and $g(n,K_3|\upp{P_3})\leq n-1$.
The required lower bound is shown by the \Lex-colored $K_n$,
 where any three vertices induce a $K_3$
 in which the two color classes are $K_2$ and $P_3$.
On the other hand, in the rainbow $K_3$ all color classes are
 single edges $K_2\notin \upp{P_3}$.
Thus, the upper bound follows from the fact $\AR(n,K_3)=n$,
 that was proved in \cite{ESS-75}.
\epf

\btm
For every\/ $n\geq 5$,\/ $\knkk=g(n,K_4|\upp{P_4})=n+1$; and\/ $g(4,K_4|\upp{2K_2})=5$,\/ $g(4,K_4|\upp{P_4})=4$.
\etm

\bpf
The smallest cases $g(4,K_4|\upp{2K_2})=5$, $g(4,K_4|\upp{P_4})=4$
 are obviously valid, as shown by a monochromatic $2K_2$ or $P_4$
 supplemented with its rainbow complementary graph using new colors.
For larger $n$, note that $2K_2\subset P_4$ implies
 $g(4,K_4|\upp{2K_2})\geq g(4,K_4|\upp{P_4})$.
Then a general construction showing the lower bound $n+1$ for
 both of them is demonstrated for every $n\geq 5$
 by the coloring of $K_n$ where a Hamiltonian cycle $C_n$ is
 rainbow-colored and its complement is monochromatic in a new color.
In the big monochromatic edge class any four vertices contain
 $P_4\in\upp{2K_2}$ as a subgraph,
 hence every copy of $K_4$ is pierced by $P_4$ in $K_n$.

In the rest of the proof, for simplicity let us abbreviate
 $\knkk$ as $g(n)$.
For the upper bound $g(n)\leq n+1$ we apply induction on $n$.
We we clearly have $g(4)\leq 5$ because the rainbow $K_4$ with
 6 colors is not pierced.
So, let $n>4$ be any larger integer, assuming that the theorem
 has been proved for $n-1$.

Consider an extremal coloring of $K_n$.
From \cite{ar-M} we recall the notion of critical colors:
 A color $c$ is critical at a vertex $v$ if all edges of color $c$
  are incident with $v$.
It means that the color class in question is either a star
 centered at $v$,
 or a single edge, which is then critical at its both ends.
The theorem follows by the simple induction $g(n)\leq g(n-1)+1$
 if there is a vertex incident with at most one critical color.

We observe that each vertex has exactly two incident critical colors.
Suppose for a contradiction that $v_0$ has at least three;
say, the edges $v_0v_i$ for $i=1,2,3$ have color $i$ critical.
In this case consider the copy of $K_4$ induced by $\{v_0,v_1,v_2,v_3\}$.
The colors $1,2,3$ do not occur inside $\{v_1,v_2,v_3\}$ because
they are critical at $v_0$, and there is not enough space
inside $\{v_1,v_2,v_3\}$ to accommodate a $2K_2$ in another color.
Hence, the coloring does not satisfy the requirements for $g$.

In fact the above argument gives more:
The two critical color classes at each $v$ are single edges,
because it does not matter if $v_0v_2$ and/or $v_0v_3$ has
color 1, still there is no room for a monochromatic $2K_2$
inside $\{v_1,v_2,v_3\}$.
So, each critical color is an edge, being critical at its both ends.
Thus, the number of critical edges equals $n$, and those edges form a graph in which each component is a cycle.

If one of those cycles is a triangle, say $C=v_1v_2v_3$,
then we select any other $v_0$ and observe that
$\{v_0,v_1,v_2,v_3\}$ violates the conditions on $g$
because the colors incident with $v_0$ are distinct from
the ones connecting $v_1,v_2,v_3$ to each other, and
their common $v_0$ does not make $2K_2$ possible.
Hence, each cycle has length at least 4.

We prove that all non-critical edges have the same color,
consequently there are just $n+1$ colors altogether.
Let $v_1,v_2,v_3$ be any three consecutive vertices along a cycle
(hence $v_1v_2$ and $v_2v_3$ are two adjacent critial edges),
and let $v_0$ be any fourth vertex.
Now the only way to satisfy the conditions for $g$
on $\{v_0,v_1,v_2,v_3\}$ is that the two edges
$v_0v_2$ and $v_1v_3$ have the same color, because the two
critical colors are not incident with $v_0$.
It means that all non-critical edges incident with any one
vertex (in the present example $v_2$) have the same color.
But then this color is the same for all vertices, because
the non-critical edges form a connected graph
(of diameter 2) whenever $n$ exceeds 4.
\epf

\btm   \label{t:uppK3}
For every\/ $n$,\/ $\knkh=n$, and\/ $g(n,K_4-e|\upp{K_3})=1$.
\etm

\bpf
We begin with the simpler case, $K_4-e$.
Of course, the monochromatic $K_n$ satisfies the conditions,
 hence $g(n,K_4-e|\upp{K_3})\geq 1$.
Consider any edge coloring $\psi$ of $K_n$ with more than one color.
Then there exist three vertices $x,y,z$ such that
 $\psi(xy)\neq \psi(yz)$.
Take a fourth vertex $w$.
The color of $wy$ differs from at least one of $\psi(xy)$ and
 $\psi(yz)$, say $\psi(wy)\neq \psi(yz)$.
Now the copy of $K_4-e$ whose edges are $wy,wz,xy,xz,yz$ contains
 no monochromatic triangle, consequently $g(n,K_4-e|\upp{K_3})< 2$.

A rainbow star $K_{1,{n-1}}$ supplemented with a monochromatic
 $K_{n-1}$ in a new color on its leaves is a piercing because
 every $K_4\subset K_n$ has at least three vertices in the
 monochromatic part.
This yields $\knkh \geq n$ for every $n\geq 4$.

For the proof of the upper bound $\knkh \leq n$, we apply
 induction on $n$.
The inequality is clearly true for $n=4$ because there can occur
 at most three further colors beside the monochromatic $K_3$
 needed for piercing.
We may also artificially extend the bound for smaller values of $n$,
 admitting $g(n,K_4|\upp{K_3})=\binom{n}{2}$ for $n=1,2,3$
 (i.e., $g(1)=0$, $g(2)=1$, $g(3)=3$).

Let $n>4$ and consider any edge coloring $\psi$ of $K_n$ with $n$
 or more colors.
Since $\AR(n,K_3)=n$, we can find a rainbow $K_3$, say with the
 three vertices $x,y,z$ where the colors are $\psi(xy) = a$,
 $\psi(yz) = b$, $\psi(zx) = c$.
Supplementing this $K_3$ with any further vertex $v$ there must
 occur a monochromatic triangle.
Hence, the other $n-3$ vertices can be partitioned into three sets
 $A,B,C$ where $A=\{v \mid \psi(vx) = \psi(vy) = a\}$ and $B,C$
 are defined analogously with respect to colors $b$ and $c$.

We assume without loss of generality that $A\neq\es$; one or both
 of $B$ and $C$ may be empty.
However, any two nonempty sets from $A,B,C$ are joined only by
 edges whose colors are from $a,b,c$.
Indeed, if $v\in A$, $w\in B$, and $\psi(vw)=d\notin\{a,b,c\}$
 held, then $vxzw$ would be a rainbow 4-cycle (having colors
 $a,c,b,d$ in this order), leading to the contradiction that the
 $K_4$ induced by these four vertices contains no monochromatic $K_3$.

The fact that no colors other than $a,b,c$ can occur on edges of
 types $A$--$B$, $A$--$C$ and $B$--$C$, has the particular
 implication that if $|A|=1$ then
  $\knkh\leq g(n-1,K_4|\upp{K_3})+1$,
 because removing the single vertex of $A$ we can destroy only one
 color, namely that of the edge from $A$ to $z$.

Consider now the set $A$, and assume $|A|>1$.
If it has a vertex $v$ such that $\psi(vz) = d$ is a new color,
 then let $w$ be any vertex in $A-v$.
We look at the $K\cong K_4$ induced by $\{x,z,v,w\}$.
Removing $v$ or $w$ from $K$, there remain the distinct colors
 $a,c$, hence $K$ is not pierced in this way.
If piercing is achieved by removing $z$, then we must have
 $\psi(vw) = a$; and if it is by the removal of $x$, then
 color $d$ remains present, hence the piercing color is then
 forced to be $\psi(vw) = d$.
Consequently, a new color $d$ between $A$ and $K$ incident with a
 vertex $v$ forces $v$ to be incident with only the colors
 $a,b,c,d$ in the entire $K_n$.
Thus, $\knkh \leq g(n-1,K_4|\upp{K_3}) +1 \leq n$ follows
 by induction applied to $K_n-v$.

Otherwise, if all edges between $A$ and $K$ are in the colors
 $a,b,c$, consider the subgraph induced by $A$.
Its coloring is supposed to be a piercing, hence at most
 $g(|A|) \leq |A|$ colors occur inside $A$, due to the
 induction hypothesis.
Moreover only the colors $a,b,c$ connect $A$ with $B\cup C\cup K$.
Applying a large induction jump now from $K_n -A$ to $K_n$,
 we obtain $\knkh \leq g(n-|A|,K_4|\upp{K_3}) + |A| \leq n$.
\epf

Replacing $K_3$ with $P_3$, the growth order jumps from linear
 to $n\sqrt{n}$; here we can derive estimates which are tight
 apart from a multiplicative constant.

\btm   \label{t:uppP3}
 $g(n,K_4|\upp{P_3}) = \Theta(n^{3/2})$, and also the smaller
 class\/ $\upp{\{K_3,P_4\}}$ generates the same growth order\/
  $g(n,K_4|\upp{\{K_3,P_4\}}) = \Theta(n^{3/2})$.
\etm

\bpf
Since $\upp{P_3}$ is obtained from $\upp{\{K_3,P_4\}}$ by
 supplementing the latter with $P_3$ and $K_{1,3}$,
 we have $g(n,K_4|\upp{P_3}) \geq g(n,K_4|\upp{\{K_3,P_4\}})$.
Then the upper bound $O(n^{3/2})$ is easily obtained by
 $g(n,K_4|\upp{P_3})
 < \Lar{K_4} \leq n\ceil{\sqrt{2n}\,}$, applying the
  corresponding estimate from \cite{ar-M}.
For the lower bound we take $C_4$-free bipartite graphs $B_n$
 whose numbers of edges grow with $n^{3/2}$.
(For instance, if $q$ is a prime power, then the point-line
 incidence bipartite graph on $2(q^2+q+1)$ vertices is
 $(q+1)$-regular and $C_4$-free.)
Let $\psi$ be the edge coloring of $K_n$ with rainbow $B_n$
 and its monochromatic complement in a new color.
If four vertices are distributed as $3+1$ or $4+0$ between the
 two classes of $B_n$, then the corresponding copy of $K_4$
 contains a monochromatic $K_3\in\upp{P_3}$.
And if the distribution is $2+2$, then the induced subgraph of
 $B_n$ cannot be $K_{2,2} = C_4$, hence an edge from the big
 color class must occur, which is extendible to a monochromatic
 $P_4\in\upp{P_3}$.
This implies $g(n,K_4|\upp{\{K_3,P_4\}}) \geq e(B_n)+1$, completing the
 proof of the theorem.
\epf

\subsection{The method of graph packing}

The method introduced in this subsection is based on the
 following notion.

\bdf
\nev{Packing graphs} \
Let\/ $G_1$ and\/ $G_2$ be two graphs on\/ $p$ vertices.
We say that\/ $G_1$ and\/ $G_2$ pack (in\/ $K_p$) if there exist
 bijections\/ $\eta_i: V(G_i)\to V(K_p)$ for\/ $i=1,2$ such that
 the induced mappings on the edge sets satisfy\/
 $\eta_1(E(G_1)) \cap \eta_2(E(G_2)) = \es$.
\edf

\paragraph{Notation.}

For a graph $F$ we denote by $\upp{F}$ the family of all graphs
 that contain $F$ as a subgraph.
Analogously, for a collection $\{F_1,\dots,F_k\}$ of graphs,
 $\upp{F_1,\dots,F_k}$ denotes the family of those graphs which
 contain at least one $F_i$ as a subgraph ($1\leq i\leq k$).

\bpn
\nev{Packing Lemma} \
Let\/ $p$ be an integer, and let\/ $H$ and\/ $Q$ be given graphs,
 such that\/ $|H|,|Q|\leq p$.
If all\/ $H$-free graphs of order\/ $p$
 (also those having some isolated vertices) pack with\/ $Q$
 in\/ $K_p$, then\/ $g(n,K_p | \upp{Q}) \geq \ex(n,H) +1$.
More generally, if\/ $\cH$ is a family of graphs, such that
 all graphs of order\/ $p$ not containing any members of\/ $\cH$ pack
  with\/ $Q$ in\/ $K_p$, then\/ $g(n, K_p | \upp{Q}) \geq \ex(n, \cH ) +1$.
\epn

\bpf
We take a rainbow extremal graph for $H$,
 with $\ex(n,H)$ edges and that many distinct colors, and
 make its complementary graph monochromatic in a new color $c$.
Consider any $K\subset K_n$, $K\cong K_p$, and let
 $F_c\subset K$ be its subgraph whose edges have color $c$,
 with $V(F_c)=V(K)$.
The complementary subgraph $K-F_c$ formed by the rainbow edges is
 $H$-free, hence it packs with $Q$ in $K_p$.
Thus, $Q\subset F_c$, and consequently the color class $F_c$ in
 $K$ is a member of $\upp{Q}$, satisfying the conditions for $g$.
The situation in the more general case of $\cH$ is similar:
 we can take a rainbow extremal graph for $\cH$ with a
 monochromatic complement in a new color $c$,
 and observe that in any copy of $K_p$
 the ranbow subgraph contains no member of $\cH$, hence
 a copy of $Q$ occurs in the $c$-colored subgraph of $K$.
\epf

We will apply the following classic result of this area.

\btm {\rm (Bollob\'as and Eldridge \cite{BE-78})} \
\label{t:BE}
Suppose that\/ $G_1$ and\/ $G_2$ are graphs with\/ $n$ vertices,\/
 $\Delta(G_1), \Delta(G_2) <n-1$,\/ $e(G_1) +e(G_2) \leq 2n-3$,
 and\/ $(G_1,G_2)$ is none of the following pairs:\/
$(2K_2,K_1 \cup K_3)$, $(2K_1 \cup K_3,K_2 \cup K_3)$,
 $(3K_2,2K_1 \cup K_4)$, $(3K_1 \cup K_3,2K_3)$,
 $(2K_2 \cup K_3,3K_1 \cup K_4)$, $(4K_1 \cup K_4,K_2 \cup 2K_3)$,
  $(5K_1 \cup K_4,3K_3)$.
Then\/ $G_1$ and\/ $G_2$ pack.
\etm

An application of this theorem will be given in
 Proposition \ref{p:sparse-H} below.

\bpn   \label{p:comp-match}
\nev{Basic lower-bound construction} \
For a given\/ $p$ and sufficiently large\/ $n$, let\/ $R$ be a\/
 $(2p-2)$-regular graph with girth\/ $g > 2p$.
Let\/ $\psi$ be the edge coloring of\/ $K_n$ composed from rainbow\/
 $R$ and its monochromatic complement assigned with a new color\/ $c$.
Then every\/ $K\subset K_n$,\/ $K\cong K_{2p}$, contains a
 monochromatic\/ $pK_2$ (in color\/ $c$).
As a consequence,\/ $g(n, K_{2p} | \upp{pK_2}) \geq (p-1)n+1$ holds
 for every\/ $n\geq n_0(p)$.
\epn

\bpf
Since the girth of $R$ exceeds $2p$, any $K$ (copy of $K_{2p}$)
 induces an acyclic graph $R_K$ with edge set $E(R)\cap E(K)$,
 that has at most $2p-1$ edges.
Further, $\Delta(K_R)\leq 2p-2$ as $K_R\subset R$.
Now we have $e(K_R) +e(pK_2) \leq 2p- 1 + p = 3p-1
 \leq 2(2p) - 3 = 4p-3$ for $p \geq 2$, and of course
 $\Delta(pK_2)=1$.
Thus, none of the forbidden cases Theorem \ref{t:BE} occurs,
 so a monochromatic $pK_2$ packs with $K_R$.
This proves the lower bound on $g(n, K_{2p} | \upp{pK_2})$
 because $(2p-2)$-regular graphs of girth larger than $2p$
 exist whenever $n$ is large enough with respect to $p$,
 by the theorem of \erd\ and Sachs \cite{ES-63}. 
\epf

The lower-bound construction is applicable not only for matchings
 but also for many more graphs.

\btm   \label{p:sparse-H}
Let\/ $H$ be a graph with\/ $p$ vertices (isolates allowed) and
 at most\/ $p-2$ edges.
If\/ $n$ is sufficiently large,\/ $n>n_0(p)$, then\/
 $g(n,K_p|\upp{H})>(p-2)n/2$.
\etm

\bpf
Similarly to the previous proof, we take a rainbow graph $R$ of
 order $n$, which is $(p-2)$-regular if $pn$ is even, or one
  vertex has degree $p-3$ and all the others have degree $p-2$
 if $pn$ is odd; and the girth of $R$ is larger than $p$.
If $n$ is large enough with respect to $p$, such a graph $R$ exists.
Make the complementary graph $\overline{R}$ monochromatic in a
 new color.
Then any $p$ vertices induce a forest in $R$, which packs with
 $H$ due to Theorem \ref{t:BE}.
Indeed, in the list of forbidden pairs of Theorem \ref{t:BE},
 precisely the seven graphs $2K_2$, $2K_1 \cup K_3$, $3K_2$,
  $3K_1 \cup K_3$, $2K_2 \cup K_3$, $4K_1 \cup K_4$,
   $5K_1 \cup K_4$ are sparse enough to play the role of $H$,
 but each of these has a pair that contains a $K_3$, hence
 cannot arise from $R$.
Consequently, in every copy of $K_p$, the subgraph induced by the
 big color class is a member of $\upp{H}$, and so the conditions
 for $g(n,K_p|\upp{H})$ are satisfied.
\epf

One may also observe that the lower-bound construction on
 $\knkk$ is a rather particular case of the above, putting $p=4$.

Next, we prove that the lower bound is tight apart from a
 muliplicative constant of 2.

\bpn
For every\/ $p$ and\/ $n$ we have\/
 $g(n, K_{2p} | \upp{pK_2}) <(2p-2)n$.
\epn

\bpf
The assertion is clearly valid for $n=2p$, because there is
 no room for $2p(2p-2)$ colors in $K_{2p}$.
For larger $n$ we prove the inequality $g(n, K_{2p} |
 \upp{pK_2}) \leq g(n-1, K_{2p} | \upp{pK_2}) + 2p-2$,
 hence applying induction.

Let $\psi$ be any coloring of $K_{2p}$ in $K_n$ with
 $g(n, K_{2p} | \upp{pK_2})$ colors.
There can be no vertex incident with $2p-1$ critical colors
 (or more), because a rainbow star
 of that many critical edges would not admit
 a $K_{2p}$ supergraph containing a monochromatic $pK_2$.
Hence, removing any one vertex we destroy at most $2p-2$ colors,
 and consequently obtain the required recursion.
\epf

\subsection{Blockers and anti-packers}

In this subsection we further elaborate the use of graph packing as a tool in investigating $g(n,G|\cF)$.

\bdf   \nev{Blocker}
For a given\/ $p \geq 3$ and a graph\/ $G$ on at most\/ $p$ vertices,
 a graph\/ $H$ is called\/ $(G,p)$-blocker if, for any placement of
 any copy of\/ $G$ and any copy of\/ $H$ in\/ $K_p$,\/
 $|E(G) \cap E(H)| \geq 2$ holds.
The graph\/ $H$ is a minimal\/ $(G, p)$-blocker if no subgraph of\/ $H$
 is a\/ $(G,p)$-blocker.
Define\/ $B(G, p)$ as the family of all minimal\/ $(G,p)$-blockers.
\edf

\bdf   \nev{Anti-packer}
For a given\/ $p \geq 3$ and a graph\/ $G$ on at most\/ $p$ vertices,
 a graph\/ $H$ is called\/ $(G,p)$-anti-packer if there is no packing
 of\/ $H$ and\/ $G$ in\/ $K_p$;
 that is, for any placement of
 any copy of\/ $G$ and any copy of\/ $H$ in\/ $K_p$,\/
 $|E(G) \cap E(H)| \geq 1$ holds.
The graph\/ $H$ is a minimal\/ $(G, p)$-anti-packer if no subgraph of\/ $H$
 is a\/ $(G,p)$-anti-packer.
Define\/ $\AP(G, p)$ as the family of all minimal\/ $(G,p)$-anti-packers.
\edf

Surprisingly, sometimes anti-packers alone are strong enough to
 determine the asymptotic behavior of $g(n,K_p|\upp{G})$.
This fact is shown in the next result.

\btm   \label{t:chi-ap}
Let\/ $G$ be a graph of order at most\/ $p$.
 \begin{itemize}
  \item[$(i)$] If there exists a bipartite anti-packer\/
   $H\in\AP(G,p)$, then\/ $g(n,K_p|\upp{G})=O(n^{2-\epsilon})$
   for some\/ $\epsilon = \epsilon(H) >0$.
  \item[$(ii)$] If the minimum chromatic number of anti-packers\/
   $H\in\AP(G,p)$ is\/ $q\geq 3$, then\/
    $g(n,K_p|\upp{G})=(1+o(1)\,\ex(n,K_q)$.
 \end{itemize}
\etm

\bpf
We begin with the upper bound, concerning which $(i)$ and $(ii)$
 can be handled essentially together.
Let $H$ be an anti-packer of minimum chromatic number, $q\geq 2$.
We are going to apply Theorem \ref{t:CCCL} and its notation,
 here writing $r:=q$.
Let $\psi$ be any edge coloring of $K_n$, using more than
 $\ex(n,K_q) + h(n)$ colors, where $h(n)=o(n^2)$, if $q\geq 3$.
If $q=2$, we assume that $\psi$ uses
 $\ex(n,K_{s,s})+1$ colors, for a sufficiently large $s$.
This means $O(n^{2-\epsilon})$ colors with $\epsilon=\epsilon(s) >0$
 if $q=2$, and a quadratic number of colors if $q\geq 3$.

Since $\chi(H)=q$, there exists a clean rainbow complete
 $q$-partite graph $Q\subset K_n$, $Q\cong K_{t,\dots,t}$,
 with vertex classes $V_1,\dots,V_q$ such that $H\subset Q$ holds.
We fix the position of $H$ in $Q$, and extend $H$ to a copy
 $K\cong K_p$ with $V(H)\subset V(K)\subset V(Q)$.
We claim that this $K$ violates the conditions on $g$, namely
 $K$ does not contain $G$ as a monochromatic subgraph.
Indeed, $G$ must share an edge $e$ with
 $H$, for otherwise $G$ would pack with $H$,
 contradicting the assumption that $H$ is an anti-packer.
Further, it is impossible for a monochromatic $G$ to have more
 than one edge in $Q$, because the edges of $Q$ have mutually
 distinct colors.
But then $G-e$ has all its edges inside the color classes $V_i$
 of $Q$, whose colors
 cannot be the same as $\psi(e)$, because $Q$ is clean.
This contradicts the assumption that $G$ is monochromatic.

For the other side of $(ii)$ we take the \rtem$(n,q-1)$ coloring
 of $K_n$, whose core is a rainbow $R\cong T_{n,q-1}$.
Recall that $R$ is a $(q-1)$-partite graph.
Consider any copy $K\subset K_n$ of $K_p$, and
 let $H$ be its intersection
 with $R$, where $V(H) = V(K)$ and $E(H) = E(K)\cap E(R)$.
Since $\chi(H)\leq \chi(R) = q-1$ and the minimum chromatic,
 number in $\AP(G,p)$ is $q$, it follows that $H$ packs with $G$
 inside $K$.
But the complement of $R$ is monochromatic in \rtem$(n,q-1)$.
Thus, every copy of $K_p$ contains a monochromatic copy of $G$,
 implying $g(n,K_p|\upp{G}) > |E(R)| = \ex(n,K_q)$.
\epf

\bcr   \label{c:conn-p/2}
If\/ $G$ is a connected graph of order at least\/ $p/2 + 1$, or
 disconnected but has a component on at least\/ $p/2 + 1$ vertices,
  then\/ $g(n,K_p|\upp{G})=o(n^2)$.
\ecr

\bpf
The complete bipartite graph $K_{\ceil{p/2},\floor{p/2}}$ is an
 anti-packer for $G$, because $p/2+1>(p+1)/2\geq\ceil{p/2}$,
 therefore even the larger vertex class is not big enough for
 accommodating $G$.
\epf

In Theorem \ref{t:chi-ap} the chromatic number was the
 relevant parameter.
In general, the \arr\ number is useful from both sides, as shown next.

\bpn   \label{p:AP-g-B}
$\AR(n, \AP(G, p)) \leq g(n,K_p|\upp{G}) < \AR(n, B(G,p))$.
\epn

\bpf
Let $\psi$ be any edge coloring of $K_n$ using at least
 $\AR(n,B(G,p))$ colors.
By definition, there is a rainbow copy $H'\subset K_n$ of some
 $H\in B(G,p)$ under~$\psi$, say inside the complete subgraph
 $K'\subset K_n$, $K'\cong K_p$.
In order to satisfy the conditions of the $g$ function, this $K'$
 should contain a monochromatic copy $G'$ of some member of $\upp{G}$.
However, since $H$ is a $(G,p)$-blocker, it blocks not only $G$
 but also all members of $\upp{G}$.
Hence, if a required $G'$ existed, we would then have
 $|E(G') \cap E(H')| \geq 2$, yielding at least two distinct
 colors in $G'$, a contradiction.
This implies $g(n,K_p|\upp{G}) < \AR(n, B(G,p))$.

To prove the lower bound $g(n,K_p|\upp{G}) \geq \AR(n, \AP(G, p))$,
 choose an edge coloring of $K_n$ with $\AR(n,\AP(G,p))-1$
 colors, without any rainbow subgraphs $H\in \AP(G,p)$.
Select just one edge from each color class, and make the rest
 monochromatic with one extra color.
In this way an edge coloring $\psi$ of $K_n$ with $\AR(n,\AP(G,p))$
 colors is obtained, in which any $K_p$ contains from the
 original colors only a subgraph that packs with $G$.
Thus, the big monochromatic class in any $K_p\subset K_n$ is a
 member of $\upp{G}$, and the requirements of function $g$ are
 satisfied.
\epf

\bex
If\/ $p=4$ and\/ $G=K_3$, then\/ $\AP(G,4)=\{2K_2\}$ and\/
 $B(G,p)=\{C_4\}$.
In this case\/ $\AR(n, \AP(G, p)) = 2$,\/ $g(n,K_p|\upp{G}) = n$,\/
 $\AR(n, B(G,p)) = \floor{4n/3}$.
\eex

This example already indicates that $\AR(n, \AP(G, p))$ and
 $\AR(n, B(G,p))$ are not very close to each other.
Nevertheless, they may serve as useful estimates on
  $g(n,K_p|\upp{G})$ from at least one side.
We begin with a case whose $g$ values are farthest from the
 \arr\ numbers of blockers.
Let us recall that $K_2\notin\cF$ is assumed in the definition of
 $g(n,G|\cF)$, therefore the notation $g(n,K_2|K_2)$ is meaningless.

\btm
Apart from the single exception that\/ $g(4,2K_2|2K_2)=3$,
 all graphs\/ $G\neq K_2$ have\/ $g(n,G|G)=1$ for all\/ $n$.
\etm

\bpf
Using more than three colors in $K_4$ would yield a 2-colored
 $2K_2$, which is not allowed by the conditions on $g$.
On the other hand, the proper edge 3-coloring of $K_4$ paints
 each $2K_2$ monochromatic.

Assume $G\neq 2K_2$, and suppose that an edge coloring $\psi$
 of $K_n$ makes all copies of $G$ monochromatic.
If $G$ is not a matching, we choose a $P\subset G$, $P\cong P_3$;
 let $E(P)=\{e',e''\}$.
Identifying an edge $e_1$ of $K_n$ with $e'$ in suitably
 chosen copies of $G$, the image $e_2$ of $e''$ must get
 $\psi(e_2)=\psi(e_1)$.
Next, identifying $e'$ with any of those $e_2$, any incident edge
 $e_3$ in the role of $e''$ must get $\psi(e_3)=\psi(e_2)$.
Thus, the entire $K_n$ is monochromatic under $\psi$.

The situation is similar if $G$ is a matching other than $2K_2$,
 or $G=2K_2$ and $n>4$.
In that case all edges disjoint from $e_1$ must get
 $\psi(e_2)=\psi(e_1)$, hence inducing a $K_{n-2}$
 entirely painted with $\psi(e_1)$.
But now $n-2>2$ holds, therefore any edge $e_3$ connecting $e_1$
 with the big monochromatic $K_{n-2}$ is disjoint from an edge
 of this $K_{n-2}$, and is extendible to a copy of $G$.
It follows that also the connecting edges must have the same
 color, consequently $\psi$ uses only one color.
\epf

Applying blockers, we can obtain a linear or subquadratic upper bound on
 $g(n,K_p|\upp{G})$ under some very natural simple assumptions.

\btm   \nev{Minimum and maximum degree}   \label{t:minmaxdeg}
Let\/ $G$ be a graph with\/ $p\geq 3$ vertices, and let\/ $n\to\infty$.
 \begin{itemize}
  \item[$(i)$] If\/ $G$ has no isolated vertices, i.e.,\/ $\delta(G)\geq 1$,
   then\/ $g(n,K_p|\upp{G}) = O(n)$.
  \item[$(ii)$] If\/ $G$ has isolated vertices, but\/
   $\Delta(G)\geq p/2$, then\/ $g(n,K_p|\upp{G}) =
    O(n^{2-\epsilon})$ for some\/ $\epsilon=\epsilon(G)>0$.
 \end{itemize}
\etm

\bpf
Let us begin with the proof of $(ii)$, as it is much simpler.
As in Corollary \ref{c:conn-p/2}, we see
 that $K_{\ceil{p/2},\floor{p/2}}$ is an anti-packer of $G$
 whenever $\Delta(G)\geq p/2$, because
 a maximum-degree vertex has too many neighbors to fit all in its
 vertex class inside $K_{\ceil{p/2},\floor{p/2}}$.
Thus, Theorem \ref{t:chi-ap} $(i)$ can be applied for the current $(ii)$.

Turning now to $(i)$, suppose first $\delta(G) \geq 2$.
Then $H = K_{1,p+1-\delta(G)}$ is a blocker for $G$, as any
 placement of $G$ and $H$ in $K_p$ gives $|E(G) \cap E(H)| \geq 2$.
Hence, by Proposition \ref{p:AP-g-B} we obtain
 $g(n,K_p|\upp{G}) < \AR(n,K_{1,p+1-\delta(G)})$,
  and the right-hand side is known to be linear in $n$.

Suppose next that $\delta(G)=1$ holds and $G$ is not a star.  
Assume that the number of colors is at least $\AR(n,K_{1,s})$,
 where $s$ is chosen to be large enough to force a homogeneously
 colored $K_p$ on the leaves of a rainbow $K_{1,s}$,
 as guaranteed by Theorem \ref{t:E--R}.
If the homogeneous $K_p$ is either rainbow or \Lex-colored,
 we are done as we can take just those $p$ vertices and no
  monochromatic copy of $G$ occurs in that copy of $K_p$.
Else, if $K_p$ is monochromatic, then the rainbow star
 $K_{1,p}\subset K_{1,s}$ attached to $K_p$ can be reduced to
 $K_{1,p-1}$ by dropping a vertex from $K_p$ and eliminating
 the shared color.
In this way the edges of $G$ incident with the vertex mapped to
 the center of the star have colors different from the colors of
 the edges not incident with the center, therefore
 no monochromatic $G$  can arise in the modified copy of $K_p$.
Since the required $s$ depends only on $p$, and $\AR(n,K_{1,s})$
 is linear in $n$ for every fixed $s$, we are done also in this case.

Finally, assume $G=K_{1,p-1}$.
In this last case $C_p$ is a blocker of $G$, hence
 $g(n,K_p|\upp{K_{1,p-1}}) < Ar(n ,C_p)$,
 again known to be bounded above by $O(n)$.
This completes the proof.
\epf

\subsection{Piercing by complete and nearly complete graphs}

Next, we list some blockers for complete subgraphs.

\bpn   \label{p:minblocker}
\
 \begin{itemize}
  \item[$(i)$] The graph\/ $H(t) := C_3 \cup tK_2$ with\/ $t \geq 1$
   is a minimal blocker in\/ $B(K_{p-t}, p)$ for\/ $p \geq 2t + 3$.
  \item[$(ii)$] The matching\/ $M(t) = tK_2$ with\/ $t \geq 3$
   is a minimal blocker in\/ $B(K_{p-t+2}, p)$ for\/ $p \geq 2t +4$.
  \item[$(iii)$] The path\/ $P_{2t +1}$ with\/ $t \geq 2$
   is a minimal blocker in\/ $B(K_{p - t +1}, p)$ for\/ $p \geq 2t +1$.
  \item[$(iv)$] The cycle\/ $C_{2t}$ with\/ $t \geq 2$
   is a minimal blocker in\/ $B(K_{p - t +1}, p)$ for\/ $p \geq 2t$.
  \item[$(v)$] Generalizing\/ $H(1)$, the cycle-plus-edge graph\/
   $C_{2t+1}\cup K_2$ with\/ $t \geq 1$
   is a minimal blocker in\/ $B(K_{p - t}, p)$ for\/ $p \geq 2t+3$.
  \item[$(vi)$] The graph\/ $tP_3$ with\/ $t \geq 2$
   is a minimal blocker in\/ $B(K_{p - t +1}, p)$ for\/ $p \geq 3t$.
 \end{itemize}
\epn

\bpf
By definition, each of the graphs claimed to be blockers have
 at most $p$ vertices, hence are
 subgraphs of $K_p$ under the given conditions.
Each of them is a blocker, because
 \begin{itemize}
  \item removing $t$ vertices from $H(t)$, the remaining graph
   contains $K_3$ or $2K_2$;
  \item removing $t-2$ vertices from $M(t)$, the remaining graph
   contains $2K_2$;
  \item since each of $P_{2t +1}$, $C_{2t}$, and $tP_3$ has
   $2t$ edges and maximum degree 2, removing $t-1$ vertices
   we can eliminate at most $2t-2$ edges and at least two remain;
  \item similarly, since $C_{2t+1} \cup K_2$ has $2t+2$ edges
   and maximum degree 2, removing $t$ vertices
   we can eliminate at most $2t$ edges and at least two remain.
 \end{itemize}
However, omitting any one edge from any of these graphs, there
 will be a set of $t$ or $t-2$ or $t-1$ vertices, respectively,
 that meets all but one of the edges.
For this reason, all these blockers are minimal.
\epf

The \arr\ numbers of these graphs have been studied in
 \cite{GR-16}, \cite{CLT-09}, \cite{Y-21+}, and \cite{S-04}, respectively.
More explicitly, we quote from \cite{GR-16} that the formula
  $$
    tn - \binom{t+1}{2} +2
  $$
 correctly determines both $\AR(n,K_3\cup tK_2)$ and
 $\AR(n,(t+1)P_3)$ for all $t\geq 1$, provided that $n>5t/2+6$
 for $K_3\cup tK_2$, and $n>5t+6$ for $(t+1)P_3$.
The results on $M(t)$ and $P_{2t+1}$ in \cite{CLT-09} and \cite{Y-21+} are similar.

We can use this approach to exactly determine the value of the
 function $g(n,p,q)$ for nearly half of the pairs $p,q$.
We formulate the result together with general lower and upper bounds.

\btm   \label{t:gnpq}
Assume\/ $p \geq 2t+3$.
 \begin{itemize}
  \item[$(i)$] For every fixed\/ $t\geq 1$,\/ $g(n,p,p-t) \leq tn-\binom{t+1}{2} +1$ holds for all sufficiently large\/ $n\geq n_0(t)$.
  \item[$(ii)$] 
 $g(n,p,p-t) \geq tn-\binom{t+1}{2} +1$ is valid
 for all\/ $n$, and also for\/ $p=2t+2$.
  \item[$(iii)$] If\/ $p\leq 2q-3$ and\/ $n>(5p-5q)/2+6$, then\/
 $g(n,p,q) = (p-q)n - \binom{p-q+1}{2} +1$.
 \end{itemize}
\etm

\bpf
The upper bound of $(i)$ follows from the quoted formula on
 $H(t)$, by applying Proposition \ref{p:AP-g-B} with $G=K_{p-q}$.

The lower bound of $(ii)$ is obtained by taking a $K_{n-t}$
 subgraph of $K_n$ monochromatic in color 1, and making its
 complement rainbow with higher colors.
Then every copy of $K_p$ has at least $p-t$ vertices in $K_{n-t}$,
 hence contains a subgraph from $\upp{K_{p-t}}$ in color 1.

Part $(iii)$ is a consequence of the previous two, using the fact
 that the formula on $\AR(n,H(t))$ is valid for all $n>5t/2+6$.
\epf

\brm
The principle behind the lower bound in Theorem \ref{t:gnpq} $(ii)$ can be
 extended from\/ $K_{p-t}$ to non-complete graphs\/ $G$ of order\/ $q=p-t$
 in a slightly stronger form, as follows.
Choose an integer\/ $s$,\/ $0\leq s\leq t$, color the edges incident
 with\/ $s$ specified vertices to obtain a rainbow\/ $K_n-K_{n-s}$,
 and use\/ $g(n,K_{p-s}|\upp{G})$ colors in the remaining\/ $K_{n-s}$.
In this way we obtain
 $$ \textstyle
   g(n,K_p|\upp{G}) \geq sn - \binom{s+1}{2} + g(n,K_{p-s}|\upp{G})
 $$
  for any\/ $s\leq p-q$.
\erm

A strong upper bound can be proved for $g(n,K_p|\upp{H})$ also in
 the cases where a non-complete graph $H$ differs from $K_p$ only
  with few edges.
For this purpose, we first prove the following auxiliary result.

\blm   \label{l:k-rain}
For any integer\/ $k\geq 3$, every edge coloring of\/ $K_n$ with
 at least\/ $k$ colors contains a rainbow subgraph\/ $Q$ with\/
 $e(Q)=k$ and\/ $|Q|\leq 2k-2$.
Moreover, the upper bound\/ $2k-2$ is tight.
\elm

\bpf
There are two ways to construct a graph with $k$ edges and
 more than $2k-2$ vertices, namely $kK_2$ and $P_3\cup (k-2)K_2$.
Given a coloring $\psi$ of $K_n$ with at least $k$ colors,
 let $Q$ be a rainbow subgraph
 with $k$ edges and as few vertices as possible.
Then $Q\neq kK_2$ because otherwise inserting
 an edge between two edges of the
 rainbow matching of size $k$ and deleting the edge of the same
 color from $kK_2$, the number of vertices would decrease.
(If the color of the inserted edge does not appear in $Q$, we can
 delete any third edge distinct from the two connected ones.)
Also, $Q\neq P_3\cup (k-2)K_2$ because in a similar way inserting
 an edge from the center of $P_3$ to an edge of $(k-2)K_2$
 and deleting the edge of the same color from $Q$,
  the number of vertices would decrease.
Thus, $|Q|\leq 2k-2$.

Tightness is shown by $K_{2k}$ where $\psi$ consists of a rainbow
 perfect matching and its monochromatic complement in a new color.
\epf

Now we are in a position to prove the following upper bound.

\btm   \label{t:small-comp}
Let\/ $t\geq 1$ and\/ $p\geq 2t+2$ be integers, and let\/
 $H$ be a graph with at most\/ $p$ vertices.
If\/ $e(K_p)-e(H)\leq t$, then\/ $g(n,K_p|\upp{H})\leq t+1$.
\etm

\bpf
We shall apply the principle that any graph $Q$ with
 $t+2$ edges and at most $p$ vertices is a blocker of
 any graph with $p$ vertices and at least
  $\binom{p}{2}-t$ edges.
To prove the theorem, suppose for a contradiction
 that there exists an edge coloring $\psi$
 of $K_n$ with at least $t+2$ colors, satisfying the requirements
 for $g(n,K_p|\upp{H})$.
On applying Lemma \ref{l:k-rain} we find a rainbow subgraph $Q$
 with $t+2$ edges and $|Q|\leq 2(2t+2)-2=2t+2$.
Choose a subgraph $K\cong K_p$ such that $Q\subset K\subset K_n$,
 and assume that a $H'\subset K$ with $H'\in\upp{H}$ is monochromatic.
There are at most $t$ edges of $K$ missing from $H'$, therefore
 $Q$ has at least two edges in $H'$.
But $Q$ is rainbow, while $H'$ is supposed be monochromatic.
This contradiction proves the theorem.
\epf

Returning to $g(K_p|\upp{K_q})$\,:
In the range $p\leq 2t+1$, or equivalently $p\geq 2q-1$,
 the behavior of $g(n,p,q)$ is substantially different.
We put the result in a very general form, to include all
 connected graphs.

\btm   \label{t:kpq-quad}
Let\/ $p\geq 2q-1$, and let\/ $H$ be any connected graph of order\/~$q$.
 Then
  $$
    \ex(n,K_{\ceil{\frac{p}{q-1}}}) + \ceil{\frac{p}{q-1}} - 1
    \leq g(n,p,q)
    \leq g(n,K_p|\upp{H})
   \leq (1+o(1))\, \ex(n,K_{\ceil{\frac{p}{q-1}}}) \,.
  $$
In particular, the same asymptotic formula\/
 $(1+o(1))\, \ex(n,K_{\ceil{p/(q-1)}})$ applies for paths and
 cycles, determining\/ $g(n,K_p|\upp{P_q})$ and\/ $g(n,K_p|\upp{C_q})$.
\etm

\bpf
We take the \rtdm\ coloring with the $k$-partite rainbow \tur\
 graph and $k$ monochromatic classes with distinct colors,
 where $k = \ceil{p/(q-1)} - 1$.
Then, since $(q-1)k<p$, among any $p$ vertices, at least $q$
 of them are from the same vertex class, hence they induce
 a monochromatic copy of a member of $\upp{K_q}$ which is a
 subfamily of $\upp{H}$.
This implies the validity of the lower bound.

The asymptotic upper bound is derived from Theorem \ref{t:CCCL}.
With reference to its notation, more than
 $\ex(n,K_{\ceil{p/(q-1)}}) + h(n)$ colors ($h(n)=o(n^2)$)
 guarantee the presence of a clean rainbow $K_{\ceil{p/(q-1)}}$-partite
 \tur\ graph with $q-1$ vertices in each of its classes.
Keeping $p$ of those $(q-1)\ceil{p/(q-1)}$ vertices, we obtain a
 $K_p$ which does not contain any monochromatic connected subgraph
 of order $q$, hence violating the conditions on the $g$ function.
\epf

Apart from the completely settled $g(n,4,3)=n$, the case of
 $p=2q-2$ is more complicated.
To handle it, we prove a result on the following \arr\ number.

\btm   \label{t:cick}
For every\/ $k\geq 2$ and\/ $n>n_0(k)$,
 $$
   kn - (k+2)(k-1)/2 \leq
    \AR(n,\{C_{2k+1}\cup K_2,C_{2k+2}\}) \leq kn \,.
 $$
\etm

\bpf
For the lower bound we take a complete monochromatic subgraph
 $K\cong K_{n-k}$ and make its complement $K_n-K$ rainbow with
 $k(n-k) + \binom{k}{2} = kn - \binom{k+1}{2}$ new colors.
If a rainbow cycle does not use an edge from the big color class,
 then the number of its vertices in $V(K_n)\smin V(K)$ is not
 smaller than that in $V(K)$, which means length at most $2k$.
If a rainbow cycle uses an edge from the big color class, then
 length $2k+1$ is possible, but such a cycle covers the entire
 $V(K_n)\smin V(K)$, therefore an additional $K_2$ would
  repeat the color of $K$.
Likewise, any cycle longer than $2k+1$ contains at least two
 edges of the monochromatic $K$, hence cannot be rainbow.

For the upper bound let $\psi$ be any edge coloring of $K_n$ with
 at least $kn$ colors, that contains no rainbow $C_{2k+1}\cup K_2$
 and no rainbow $C_{2k+2}$.
We first observe that $kn$ exceeds
 $\AR(n,C_{2k+1}) = ( \frac{2k-1}{2} + \frac{1}{2k} ) \, n + O(1)$
 if $n$ is large.
Thus, a rainbow cycle $C\cong C_{2k+1}$ surely is present.
Denote $X:=V(C)$ and $Y:=V(K_n)\smin V(C)$.
Inside $Y$ every edge has its color from
 $\psi(C)$, otherwise a rainbow $C_{2k+1}\cup K_2$ would exist.

Consider any vertex $y\in Y$.
At most $k$ colors different from those in $\psi(C)$ can occur
 on the edges connecting $y$ with $X$, otherwise there would be
 two of them to two consecutive vertices of $C$, and so $C$
 would be extendable to a rainbow $C_{2k+2}$, a contradiction.
Consequently the number of colors used by $\psi$ is at most
 $
   k\cdot |Y| + \binom{|X|}{2} = k(n-2k-1) + \binom{2k+1}{2}
    = kn
 $.
\epf

This result is strong enough for the determination of
 $g(n,2q-2,q)$ with high precision.

\btm   \label{t:2q-2}
If\/ $p=2q-2$ and\/ $n>n_0(p)$, then
 $$(p-q)n-\binom{p-q+1}{2} +1 \leq g(n,p,q) \leq
  \AR(n,\{C_{2p-2q+1}\cup K_2,C_{2p-2q+2}\}) \leq (p-q)n \,.$$
\etm

\bpf
For simplification, let us denote $t=p-q$.
Then the lower bound is a consequence of Theorem \ref{t:gnpq}
 $(ii)$, which is valid for all~$n$.
For the upper bound we observe by Proposition \ref{p:minblocker}
 $(iv)$--$(v)$ that both $C_{2t+1}\cup K_2$ and $C_{2t+2}$ are
 minimal blockers for $K_q$.
Hence, Proposition \ref{p:AP-g-B} and Proposition
 \ref{p:monoton} $(iii)$ imply $g(n,p,q) \leq \AR(n,B(K_q,p)) \leq
  \AR(n,\{C_{2t+1}\cup K_2,C_{2t+2}\})$.
The explicit upper bound is then obtained from Theorem \ref{t:cick}.
\epf

The obtained growth orders concerning $g(n,p,q) = g(n,K_p|\upp{K_q})$
 can be summarized as follows.

\bcr   \nev{Dichotomy}
The function\/ $g(n,p,q)$ is linear in\/ $n$ if\/ $p\leq 2q-2$, and
 quadratic otherwise.
\ecr

Using the principle of family complementaion as given in
 Proposition \ref{p:compl}, tight results can be obtained also
 for the function $f$.
We present two results of this kind as a closing of the current
 subsection.

\btm   \label{t:no-k-p-t}
Let\/ $p \geq 2t+3$ (\/$t \geq 1$),
 and let\/ $\cF$ be the family of all graphs not containing\/ $K_{p-t}$.
Then, for\/ $n>n_0(p)$ sufficiently large,\/
 $f(n,K_p|\cF) = g(n,p,p-t) +1= t(n-t) + \binom{t}{2} + 2$.
\etm

\bpf
The complement $\overline{\cF}$ of $\cF$ is the family
 $\upp{K_{p-t}}$, and of course $K_p\notin \cF$.
Thus, the formula follows by Proposition~\ref{p:compl} $(ii)$ and
 Theorem \ref{t:gnpq}.
\epf

A further interesting application of the complementarity
 principle between the functions $f$ and $g$
 is the problem considered in Theorem 5.10 of \cite{ar-G}.
Before stating the new formula, let us quote that result first.
Recall the definition $\chiF(G) := \min_{F\in\cF} \chi(G-F)$.
Let now $\trif$ denote the family of all triangle-free graphs.
Concerning the behavior of the $f$ function with respect to
 $\trif$ it was proved in \cite{ar-G} that
 \begin{eqnarray}
f(n, K_p | \trif) & = &
  (1+o(1)) \, \ex(n, K_{\chi_{_\trif}(K_p)})  \label{eq:trif} \\
  & \geq &
(1+o_n(1)) \, \ex(n, K_{(1/2 +o_p(1)) \sqrt {p \log p} }).
   \nonumber
 \end{eqnarray}
The rightmost expression was derived from the currently known
 strongest upper bound on the chromatic number of triangle-free
  graphs of given order, which is very hard to determine exactly.
We now prove the following simple formula.

\btm   \label{t:f-g-k3}
For the family\/ $\trif$ of\/ $K_3$-free graphs we have\/
$f(n, K_p | \trif) = g(n,K_p | \upp{K_3}) + 1 =
 (1+o(1))\,\ex(n,K_{\ceil{p/2}})$ for every\/ $p\geq 5$.
\etm

\bpf[First proof.]
Since the hereditary family $\trif$ contains all triangle-free
 graphs, and $K_3\notin \trif$, the complement
  $\overline{\trif}$ is the upward-closed family
 $\upp{K_3}$ of all graphs containing a triangle.
Thus, according to Proposition \ref{p:compl} we have
 $f(n, K_p | \trif) = g(n,K_p | \upp{K_3}) + 1$.
In this way asymptotic equality is implied by Theorem
 \ref{t:kpq-quad} above, putting $q=3$.

\msk

\nin
\emph{Second proof.}
Despite that triangle-free graphs can have arbitrarily large
 chromatic number, we prove that
 \begin{equation}
  \chi_{_\trif}(K_p) = \ceil{p/2}
    \label{eq:chi-trif}
 \end{equation}
  holds for every $p$.
Let $G$ be any graph of order $p$ such that $\overline{G}$ is
 triangle-free.
Suppose for a condtradiction that $\chi(G)=k<p/2$.
Let $V_1,\dots,V_k$ be the color classes in a proper vertex
 coloring with the minimum number of colors.
Since $2k<p$, some $V_i$ has size greater than 3.
Thus, $K_3\subset K_{|V_i|}\subset\overline{G}$, a contradiction.
Consequently, Equality (\ref{eq:trif}) implies the theorem
 via (\ref{eq:chi-trif}).
\epf

In the same way, the following more general result can also be proved.

\btm   \label{t:p>2q}
If\/ $p\geq 2q-1$ and\/ $q\geq 3$, then
 the family of\/ $K_q$-free graphs satisfies\/
$f(n, K_p | K_q\mathrm{-free}) = g(n,K_p | \upp{K_q}) + 1 =
 (1+o(1))\,\ex(n,K_{\ceil{p/(q-1)}})$.
\etm

\bpf
The family of $K_q$-free graphs contains $K_2$ but does not
 contain $K_p$ for any $p\geq q$.
Moreover, it is the complementary family of $\upp{K_q}$.
Thus, we can apply Theorem \ref{t:kpq-quad} together with
 Proposition \ref{p:compl} for all $p$ and $q$ satisfying the
 assumptions.
\epf

\subsection{Piercing by cycles}

We have already seen by the smallest nontrivial example $H=P_3$,
 namely by the sequence $g(n,K_3|\upp{P_3})=n-1$,
 $g(n,K_4|\upp{P_3})=\Theta(n^{3/2})$,
 $g(n,K_p|\upp{P_3})=\Theta(n^2)$ for $p\geq 5$,
  that the growth-order
 dichotomy of piercing $K_p$ by complete subgraphs does
 not remain valid in general for $g(n,K_4|\upp{H})$.
In this section we take a closer look at this situation.

We begin with a little Ramsey-type lemma and its consequence on
 cycle--tree packing.

\blm
Let\/ $a\geq b\geq 2$ be integers.
In every red-blue edge coloring of\/ $K_{a,b}$ there is a
 red\/ $2K_2$ or a blue\/ $K_{1,a}$.
\elm

\bpf
If any two edges of a graph meet, then they form a star or a
 triangle.
Hence, if the 2-coloring of $K_{a,b}$ does not contain a red $2K_2$,
 then the blue subgraph is obtained by deleting a star from it.
Consequently, at least $b-1\geq 1$ vertices have blue degree $a$.
\epf

In the present context the following implication of this lemma
 will be important, because the condition $e(G_1)+e(G_2)\leq 2p-3$
 of Theorem \ref{t:BE} is not satisfied in the case under consideration.

\bpn   \label{p:t-p-pack}
Let\/ $T$ be a tree or forest of order\/ $p\geq 3$.
If\/ $\Delta(T)\leq p-3$, then\/ $T$ and\/ $C_p$ pack in\/ $K_p$.
\epn

\bpf
If $p=3$, then $T=3K_1$, and if $p=4$, then the largest
 possible $T$ is $2K_2$.
The complements of these graphs are $C_3$ and $C_4$,
 respectively, verifying the assertion for such small $p$.

If $p\geq 5$, we may restrict attention to trees, because if $T$
 has more than one component, we may connect two leaves of
 distinct components by inserting a new edge,
 without violating the maximum-degree condition.
Viewing $T$ as a bipartite graph, let $A$ and $B$ be its vertex
 classes; say, $|A|=a\geq |B|=b$.
We have $b\geq 2$ because $T\neq K_{1,p-1}$.

Assume $T$ embedded in $K_{a,b}$, all inside $K_p$.
Color the edges of $T$ blue, and the edges of $K_{a,b}$ red.
If there exist two red edges $e,e'$, supplement them with $e$--$e'$
 paths $P_a$ and $P_b$ of $K_p$ inside $A$ and $B$, respectively.
In this way a $C_p$ is obtained, which is edge-disjoint from $T$.
Hence, it suffices to prove that such red $e,e'$ exist.

If there is no blue $K_{1,a}$, then there are red $e,e'$ by the
 preceding lemma.
On the other hand, if there is a blue $K_{1,a}$, then we have
 $a\leq \Delta(T)\leq p-3$, hence $b\geq 3$.
More than one blue $K_{1,a}$ centered in $B$ would imply the
 contradiction $K_{a,2}\subset T$.
Thus, omitting the vertex of degree $a$ from $B$, we obtain a
 subforest of $T$ with more than one vertex in each vertex class,
 without blue $K_{1,a}$, consequently containing a red $2K_2$.
This completes the proof.
\epf

\btm
Piercing by the\/ $4$-cycle for large enough\/ $n$ satisfies the
  following relations: 
 \begin{itemize}
  \item[$(i)$] $g(n,K_p|\upp{C_4}) =
   (1+o(1))\,\ex(n,K_{\ceil{p/3}})$ for\/ $p\geq 7$;
  \item[$(ii)$] $g(n,K_p|\upp{C_4}) = \Theta(n^{3/2})$ for\/ $p=5,6$;
  \item[$(iii)$] $g(n,K_4|\upp{C_4}) =
   \AR(n,K_{1,3}) - 1 = \floor{n/2} + 1$.
 \end{itemize}
\etm

\bpf
Assume throughout the proof that $n$ is sufficiently large
 with respect to $p$.

\msk

\nin
$(i)$\quad 
We know that the $g$ function for $C_4$ is not smaller than
 that for $K_4$.
Hence we only need to prove $g(n,K_p|\upp{C_4})
 \leq (1+o(1))\,\ex(n,K_{\ceil{p/3}})$.
Using more than $\ex(n,K_{\ceil{p/3}}) + h(n)$ colors, where
 $h(n)=o(n^2)$, by Theorem \ref{t:CCCL} we can find a clean
 rainbow $K\cong K_{3,\dots,3}$ with $\ceil{p/3}$ vertex classes.
Inside $K$, none of the $K_p$ subgraphs contains a monochromatic
 $C_4$, thus no coloring with so many colors can satisfy the
 conditions on $g$.

\msk

\nin
$(ii)$\quad 
The upper bound follows by the fact $\AR(n,K_{3,3}) = \Theta(n^{3/2})$,
 because extending a rainbow $K_{3,3}$ to any coloring of $K_6$,
 at most $2K_3$ plus an edge can occur as monochromatic.
For a lower bound growing with $n^{3/2}$, we take
 a rainbow \emph{bipartite} extremal graph $B=B_n$ for $C_4$,
 which is known to have at least $\frac{1}{2}n^{3/2}-O(n)$ edges.
Make the edges of the complement $\overline{B}$ monochromatic
 in a new color.
We prove that every $K\cong K_5$ subgraph of this $K_n$ has a
 $C_4$ in the large color class.

Let $X,Y$ be the vertex classes of $B$.
The vertex distribution of $K$ between them is either
 $5+0$ or $4+1$ or $3+2$.
The first two obviously contain a monochromatic $K_4$, hence
 also a $C_4$.
For the case of $3+2$, denote by $u,v,z$ and $x,y$ the
 vertices of $K$ in the two classes of $B$.
Temporarily restricting attention to the edges of $B\cap K$, now
 $x$ and $y$ have either just one common neighbor $z$, or
 no common neighbor.
There can be at most two edges from $\{u,v\}$ to $\{x,y\}$
 (one from $u$ and one from $v$).
If both $ux,vx$ are edges, then $uy,vy$ are non-edges, so that
 $yuzv$ is a monochromatic $C_4$ in $K_n-E(B)$.
The situation is similar if both $uy,vy$ are edges.
Finally, if $ux$ and $vy$ are edges of $B$, then $uyxv$ is a
 monochromatic $C_4$.

\msk

\nin
$(iii)$\quad 
The $K_4$ spanned by a rainbow $K_{1,3}$ cannot contain a
 monochromatic $C_4$; this implies that $\AR(n,K_{1,3})-1$ is
 an upper bound on $g(n,K_4|\upp{C_4})$.
The simple construction yielding the same value as lower bound
 is the rainbow matching with $\floor{n/2}$ edges and its
 monochromatic complement in a new color.
Clearly, any $K_4$ subgraph contains a $C_4$
 in the big color class.
\epf

Concerning longer cycles $C_p$, the following result determines
 the exact value for all large $n$.

\btm   \label{t:kpcp}
For every fixed\/ $p\geq 5$ we have\/ $g(n,K_p|\upp{C_p}) \leq
   \left\lfloor \frac{p-3}{2} \, n \right\rfloor + 1$
  for all\/ $n\geq 2p-5$.
Moreover, the equality\/ $g(n,K_p|\upp{C_p}) =
   \left\lfloor \frac{p-3}{2} \, n \right\rfloor + 1$
    is valid for every sufficiently large\/ $n>n_0(p)$.
\etm

\bpf
Observe that the star $K_{1,p-1}$ is a minimal blocker for $C_p$
 in $K_p$, therefore $g(n,K_p|\upp{C_p}) \leq
  \AR(n,K_{1,{p-1}})-1$ holds by Proposition \ref{p:AP-g-B}.
Jiang \cite{J-02} and Montellano-Ballesteros \cite{MB-06} proved
 the formula $\AR(n,K_{1,k}) =
 \left\lfloor \frac{(k-2)n}{2} \right\rfloor +
  \left\lfloor \frac{n}{n-k+2} \right\rfloor + \epsilon$,
 where $\epsilon\in\{0,1\}$.
If $n/2>k-2$, then the second term is equal to 1.
Thus, plugging in $k=p-1$, for all $n\geq 2p-5$ we obtain that
 $\left\lfloor \frac{p-3}{2} \, n \right\rfloor + 1$
 is an upper bound.

For the matching lower bound we can take a rainbow graph
 $G_n\subset K_n$ of girth bigger than $p$,
 with $e(G_n)=\left\lfloor \frac{p-3}{2} \, n \right\rfloor$ and
 $\Delta(G_n)=p-3$, and assign a new color $c$ to its complement.
With reference to \cite{ES-63} we know that such graphs exist
 for every $n> n_0(p)$ if $n$ is even or $p$ is odd; and a slight
 modification of the proof implies the same for odd $(p-1)n$ also.

Any set $X\subset V(G_n) = V(K_n)$ of $|X|=p$ vertices induces a
 forest, say $T$, of maximum degree at most $p-3$ in $G_n$.
According to Proposition \ref{p:t-p-pack}, $T$ and $C_p$ pack
 inside $X$, which means that a monochromatic $C_p$ in color $c$
  exists in every copy of $K_p$.
Hence, the number
 $\left\lfloor \frac{p-3}{2} \, n \right\rfloor + 1$ of colors is
  a lower bound on $g$.
\epf

This result determines $g(n,K_p|\upp{C_q})$ for $p=q$, whereas
 the range $p\geq 2q-1$ was settled asymptotically in
 Theorem \ref{t:kpq-quad}, showing that it is essentially the
 same as $g(n,K_p|\upp{K_q})$ and $g(n,K_p|\upp{P_q})$.
For the remaining range $q<p\leq 2q-2$ we prove relation with the
 \tur\ numbers $\ex(n,K_{s,t})$ of complete bipartite graphs,
  where $p-q<s\leq t<q$.
The well known estimates are $\Omega(n^{2-\frac{2}{s-1}}) \leq
 \ex(n,K_{s,t}) \leq O(n^{2-\frac{1}{s}})$, and tightness of the
 upper bound is proved for some small values of $s$.
In this way, improvements in the \tur-type problems on bipartite
 graphs may also improve the numerical bounds in the next result.

\btm   \label{t:kpcqpq}
Let\/ $p,q\geq 4$ be integers, with\/ $q<p\leq 2q-2$.
 \begin{itemize}
  \item[$(i)$] For cycles, there exist constants\/ $c_1=c_1(p,q)>0$
   and\/ $c_2=c_2(p,q)>0$ such that\/
    $c_1 n^{2-\frac{2}{p-q+1}} \leq g(n,K_p|\upp{C_q}) <
     \AR(n,K_{p-q+1,q-1}) \leq c_2 n^{2-\frac{1}{p-q+1}}$.
  \item[$(ii)$] For paths, there exists a constant\/ $c=c(p,q)>0$
   such that\/
    $g(n,K_p|\upp{P_q}) \geq c n^{2-\frac{2}{p-q+2}}$.
 \end{itemize}
\etm

\bpf
We relate the $g$-values in question to \tur\ numbers of complete
 bipartite graphs.
First, observe that $K_{p-q+1,q-1}$ is a blocker for $C_q$ in $K_p$
 in the given range of $p$ with respect to $q$, because
 $p-q+1\leq q-1$ holds, hence any copy of $C_q$ inside
 $V(K_{p-q+1,q-1})$ contains vertices in both partite sets of
 $K_{p-q+1,q-1}$.
Thus, the upper bound given in $(i)$ follows by
 Proposition \ref{p:AP-g-B}.

For the lower bound in $(i)$ we consider the family
 $$
   \cK = \{ K_{p-q,q-2} , K_{p-q+1,q-3} , \dots,
     K_{\floor{p/2}-1,\ceil{p/2}-1} \} \,.
 $$
We claim that $g(n,K_p|\upp{C_q}) > \frac{1}{2}\,\ex(n,\cK)$.
From this, the standard lower-bound proof on $K_{s,t}$ yields
 the stated growth order from below.
As a matter of fact, instead of $\frac{1}{2}\,\ex(n,\cK)$ one can also take
 the bipartite \tur\ number of $\cK$, i.e., the maxumum number
 of edges in a subgraph $H\subset K_{\floor{n/2},\ceil{n/2}}$
 such that $H$ has no $K\subset H$ with $K\in \cK$.

So, let $Q$ be either an extremal bipartite $\cK$-free graph or
 a largest edge-cut of an extremal $\cK$-free graph of order $n$.
Consider the edge coloring $\psi$ of $K_n$ where $Q$ is rainbow
 and its complement is monochromatic in a new color.
Then the number of colors is $e(Q)+1 > \frac{1}{2}\,\ex(n,\cK)$.
Let us denote by $A$ and $B$ the vertex classes of $Q$.

Consider any $K\cong K_p$ in $K_n$.
Assume that $K$ has $a$ vertices in $A$ and $b$ vertices in $B$,
 with $a+b=p$ and $a\leq b$.
(To achieve this, rename $A\leftrightarrow B$ if necessary.)
If $b\geq q$, then a monochromatic $K_q$ is present in $K$,
 hence also a monochromatic $C_q$.
Otherwise $b\leq q-1$ and $a\geq p-q+1$.
Since $Q$ is $\cK$-free, the $K_{a,b}$ induced by $K$ contains
 an edge, say $xy$ with $x\in A$, $y\in B$, such that
 $xy\notin E(Q)$.
Further, restricting $K$ to $K-x-y$, the corresponding $K_{a-1,b-1}$
 is also a member of $\cK$, hence it has an edge, say $wz$ with
  $w\in A$, $z\in B$, such that $wz\notin E(Q)$.
Hence $\psi(xy)=\psi(wz)$ is the color of the big monochromatic class.
Choose now $a'$ and $b'$ such that $2\leq a'\leq a$,
 $2\leq b'\leq b$, and $a'+b'=q$.
Taking a path $P_A$ of order $a'$ from $x$ to $w$ inside $A$, and
 a $P_B$ of order $b'$ from $z$ to $y$ inside $B$, we obtain a
  monochromatic copy of $C_q$.
Since this can be done in every copy of $K_p$, the conditions
 for function $g$ are satisfied.

The proof of $(ii)$ is essentially the same, except that now we
 can build the construction upon the more advantageous family
 $$
   \cK' = \{ K_{p-q+1,q-1} , K_{p-q+2,q-2} , \dots,
     K_{\floor{p/2},\ceil{p/2}} \} \,.
 $$
We take a rainbow bipartite graph $Q'$ which is either a
 largest bipartite cut of a $\cK'$-extremal graph, or is a
 bipartite-extremal graph for $\cK'$, of order $n$; and
  supplement it with
 monochromatic complement of a new color in $K_n$.
Let $A',B'$ be the vertex classes of $Q'$.
Since $Q'$ is $\cK'$-free, every copy $K'$ of $K_p$ either has
 at least $q$ vertices in $A'$ or in $B'$, or else contains an
 edge $xy\notin E(Q')$, $x\in A'$, $y\in B'$.
This edge extends to a monochromatic $P_p$ inside $K'$, hence
 also a $P_q$.
\epf

\brm
The proof of Theorem \ref{t:kpcqpq} provides a method by which
 fairly good asymptotic estimates can be obtained for piercing
 of\/ $K_p$ by graphs of low edge connectivity.
For instance, if\/ $G$ is\/ $2$-edge-connected but not\/ $3$-edge-connected,
 then also\/ $g(n,K_p|\upp{G})$ satisfies the inequalities
 proved above for\/ $g(n,K_p|\upp{C_q})$.
We do not develop further the theory in this direction here, we leave
 the elaboration of these principles as a subject for future research.
\erm

We conclude this subsection with a further consequence of piercing combined with complementarity.
Similarly to Theorems \ref{t:no-k-p-t} and \ref{t:f-g-k3}, also
 here we can derive a tight result for function $f$, now
  concerning cycles.


\btm
Let\/ $p\geq 5$ be an integer, and let\/ $\cF$ be the family of
 graphs which do not contain\/ $C_p$ as a subgraph.
If\/ $n>n_0(p)$ is sufficiently large with respect to\/ $p$, then\/
 $f(n,K_p|\cF) = \left\lfloor \frac{p-3}{2} \, n \right\rfloor + 2$.
As a consequence, the identity\/
 $f(n,K_p|C_p\mathrm{-free}) = \AR(n,K_{1,{p-1}})$
 is valid for large\/ $n$.
\etm

\bpf
Observe that our current $\cF$ is the complementary family of
 $\upp{C_p}$ (and vice versa).
Further, $K_2\in\cF$ and $C_p\notin\cF$.
Hence, we can apply Proposition \ref{p:compl} $(ii)$ for $\cF$,
 and the theorem follows.
\epf

\subsection{Piercing by stars}

As stars constitute another fundamental class of graphs, here we
 consider the piercing problem for the family $\upp{K_{1,q}}$
  with $q\geq 2$.
Let us recall from Theorem~\ref{t:connbip} that
 $g(n, K_p|\upp{K_{1,q}})$ is subquadratic whenever $p\leq 2q$.
Below we give more precise estimates for the entire range of $p>q$.

\btm
The\/ $g$ function for the family\/ $\upp{K_{1,q}}$ admits the following estimates.
 \begin{itemize}
  \item[$(i)$] If\/ $p=q+1$, then\/ $g(n, K_p|\upp{K_{1,q}})$ is
   linear in\/ $n$\,:
   $$
     n + \binom{p-2}{2} - 1 \leq g(n, K_p|\upp{K_{1,p-1}})
      < \AR(n,C_p) = ( \frac{p-2}{2} + \frac{1}{p-1} ) \, n + O(1) \,.
   $$
  \item[$(ii)$] If\/ $q+2\leq p\leq 2q-2$, then\/ $g(n, K_p|\upp{K_{1,q}}) < \AR(n,K_{q-1,p-q+1}) = O(n^{2-\frac{1}{p-q+1}})$,
 subquadratic in\/ $n$.
  \item[$(iii)$] If\/ $2q-1\leq p\leq 2q$, then\/ $g(n, K_p|\upp{K_{1,q}}) = O(n^{2-\epsilon})$ where\/ $\epsilon=\epsilon(p)>0$ is constant, the function\/ $g$ being
  subquadratic in\/ $n$.
  \item[$(iv)$] If\/ $p\geq 2q+1$, then\/ $g(n, K_p|\upp{K_{1,q}}) = (1+o(1))\,\ex(n,K_{\ceil{p/q}})$,
   quadratic in\/~$n$.
 \end{itemize}
\etm

\bpf

$(i)$\quad
The lower bound is obtained by modifying the \Lex\ coloring.
Inside $\{v_1,\dots,v_{p-1}\}$ we take a rainbow $K_{p-1}$, and
 for every $v_i$ with $p\leq i\leq n$ we define a private color
 assigned to all edges connecting $v_i$ with vertices of smaller index.

The upper bound follows by Proposition \ref{p:AP-g-B}, observing
 that $C_p$ is a blocker for $K_{1,p-1}$ in $K_p$.

\msk

\nin
$(ii)$\quad
Also here, we apply Proposition \ref{p:AP-g-B}.
Generalizing the above example of $C_p$ that was taken for
 $p-q=1$, observe that every spanning subgraph $H\subset K_p$
  with $\delta(H)=p-q+1$ is a blocker of $K_{1,q}$.
Hence, we can take $H=K_{q-1,p-q+1}$ whenever $p-q+1\leq q-1$,
 which specifies exactly the current range $p\leq 2q-2$.

\msk

\nin
$(iii)$\quad
Here we apply Theorem \ref{t:RCL} with the parameters $b=q$
 and $r=2$.
Letting then $a=a(q,2)$, the theorem guarantees that if an
 edge coloring $\psi$ of $K_n$ uses more than $\ex(n,K_{a,a})
 = O(n^{2-1/a})$ colors, then a clean rainbow $K\cong K_{q,q}$ occurs.
Since the colors used inside the vertex classes do not occur on
 the edges of $K$, a copy of $K_p$ inside $V(K)$ does not give
 enough room for a monochromatic $K_{1,q}$.

\msk

\nin
$(iv)$\quad
Lower bound is obtained from the $\rtdm(n,\ceil{p/q}-1)$
 coloring pattern.
In any copy of $K_p$, at least one vertex class contains at least
 $q+1$ vertices, hence not just $K_{1,q}$ but even $K_{q+1}$ is
 present as a monochromatic subgraph.
From the other side, Theorem \ref{t:CCCL} implies that in every
 edge coloring of $K_n$
 with more than $\ex(n,K_{\ceil{p/q}})+h(n)$ colors
 contains a clean rainbow $K\cong K_{q,\dots,q}$ with $\ceil{p/q}$
 vertex classes.
Since $q\cdot\ceil{p/q} \geq p$, inside $V(K)$ a copy of $K_p$
 violating the consitions of $g$ is found.
\epf

\section{Concluding remarks and open problems}
\label{s:conclude}

This present work is the third and last part of our trilogy on
 some generalizations of \arr\ theory.
The starting point was that in a rainbow edge-colored graph
 every color class is a single edge.
\Arr\ theorems ask for a tight lower bound on the number of
 colors that guarantees, in every allowed coloring, the presence
 of a subgraph of given type in which every color class consists
 of just one edge.
Our general approach is to generalize the ``single edge''
 color class to specified conditions on the color classes
 of the required subgraph.

The first part of the trilogy, \cite{ar-G}, introduced the
 general problem and developed the theory for conditions
 expressible in terms of hereditary graph properties.
The second part, \cite{ar-M}, defined seven natural classes of
 subproblems, some of them expressible inside the ``graph classes''
 scenario, and some others related to popular types of graph coloring.
This current third part extends the panorama to graph properties
 which are not supposed to be hereditary anymore.
We developed some useful methods for this general approach, and
 established asymptotically tight results in several cases.
We also introduced a new \arr-type invariant, showed that it is
 strongly related to the former one by family complementation,
 and derived many further tight asymptotic results.

Almost all formulas derived here deal with the functions
 $f(n,H|\cF)$ and  $g(n,H|\cF)$ in the case of $H=K_p$, the
 complete graph of given order.
In the first and second parts of the trilogy, many results are
 given for $f(n,H|\cF)$ with non-complete $H$, but mainly with
 particularly defined $\cF$ (a considerable number of them, though).
A very natural direction of future research opens:

\bpm
\
\begin{itemize}
 \item[$(i)$] Develop the general theory to the involvement of graphs\/ $H$ other than complete graphs.
 \item[$(ii)$] Develop the theory to the involvement of families\/ $\cH$ of graphs richer than a single graph\/ $H$.
\end{itemize}
\epm

Moreover, so far the edge colorings of complete graphs have been considered.
However, \tur-type extremal theory has already been investigated
 in the literature for the restricted universe of bipartite graphs
 and also for the generalizations to hypergraphs, and more.
These directions offer a wide field of potential research in
 \arr\ theory, too.

\bpm
Develop the theory for
 \begin{itemize}
  \item bipartite graphs;
  \item directed graphs;
  \item signed graphs;
  \item uniform hypergraphs.
 \end{itemize}
\epm

\subsection{General problems on growth order and function relations}

Let us begin this subsection with some simple questions
 concerning the behavior of function $g$.

\bpm
Determine a condition which is necessary and sufficient for\/
 $g(n,G|\cF)>0$.
  (Cf.\ Remark \ref{r:g-pos}.)
\epm

\bpm
Characterize the pairs\/ $(G,\cF)$ for which\/ $g(n,G|\cF)=1$,
 or a little further\/ $g(n,G|\cF)=k$ for\/ $k$ small.
\epm

\bpm
Determine a condition which is necessary and sufficient for\/
 $g(n,G|\cF) = f(n,G|\overline{\cF})-1$.
  (Cf.\ Proposition \ref{p:compl}.)
\epm

It is worth emphasizing that the tight relation between
 $g(n,G|\cF)$ and $f(n,G|\overline{\cF})$ allows translation of
 results from $f$ to $g$ and also the other way round.
Some examples of this kind have been given in Section \ref{s:pierc},
 but many more can be derived using this tool.

For the further discussion let us recall that $f(n,G|\cF)$ is
 defined if and only if $K_2\in \cF$,
 and $g(n,G|\cF)$ is defined if and only if $K_2\notin \cF$.
So, given any $\cF$, assuming that $\phi(n,G|\cF)$ is defined
 for $\phi\in\{f,g\}$, it is uniquely determined whether
 $\phi=f$ or $\phi=g$.
In the sequel $\phi$ stands for any of $f$ and $g$, provided that
 $\phi(n,G|\cF)$ is defined.

The most general problem concerning $\phi(n,G|\cF)$
 is the following.

\bpm
Determine the families\/ $\Phi_f$ and\/ $\Phi_g$ of functions\/
 $\eta=\eta(n)$ for which there exists a graph\/ $G$
 and a family\/ $\cF$ of graphs such that
 $$
   \lim_{n\nti}\phi(n,G|\cF)/\eta(n) = 1
 $$
  as\/ $n\nti$.
\epm

Of course, the trivial identically zero function $P_0$
 (i.e., $P_0(n)=0$ for all $n$)
 is a member of $\Phi_g\smin\Phi_f$.
But we do not know which other functions are in the symmetric
 difference of $\Phi_f$ and $\Phi_g$.

\bpm   \label{pr:fi-diff}
Determine\/ $\Phi_f\smin\Phi_g$ and\/ $\Phi_g\smin\Phi_f$.
\epm

From now on we assume without further mention that $\phi$
 determines a function different from $P_0$.

It is not known whether every $\phi(n,G|\cF)$ is proportional to
 a real power of $n$ at least in the asymptotic sense.

\bpm
Prove or disprove that for every graph\/ $G$ and family\/ $\cF$ of
 graphs there exist constants\/ $c=c(G,\cF)$ and\/ $d=d(G,\cF)$
 such that
  \begin{equation}   \label{eq:limit}
   \lim_{n\to\infty} \frac{\phi(n,G|\cF)}{c\cdot n^d} = 1 \,.
  \end{equation}
\epm

Motivated by this problem, we introduce the following notion.

\bdf   \nev{Feasible pair}
We call a pair\/ $(c,d)$ of real numbers\/ $c>0$ and\/ $0\leq d\leq 2$
  feasible if there exists a graph\/ $G$ and a family\/ $\cF$ of
 graphs for which\/ $(\ref{eq:limit})$ holds.
In that case we say that\/ $G$ and\/ $\cF$ determine\/ $(c,d)$.
\edf

More directly, we raise:

\bcj
If\/ $\cF$ is finite, and\/ $\phi(n,G|\cF)\to\infty$, then\/ $G$ and\/
 $\cF$ determine a feasible pair\/ $(c,d)$ for every\/ $G$.
\ecj

Currently we expect that the assumption $|\cF|<\infty$ is
 superfluous in the conjecture.

\bpm
Characterize the feasible pairs\/ $(c,d)$ for\/ $\phi\in\{f,g\}$.
\epm

Translating to this language the results available in the \arr\
 literature, many examples of feasible pairs can be obtained.
In particular, since $\AR(n,G)$ is expressible in terms of
 function $f$, and under some conditions the former has the
 same asymptotics as $\ex(n,\cF')$ with a suitably constructed
 $\cF'$, every \tur-type theorem exhibiting an explicit asymptotic
 formula provides us with a feasible pair.

In the subquadratic world---which means $d'<2$ in terms of
 feasible pairs $(c,d)$, or bipartite forbidden graphs in \tur-type
 extremal theory---the following problem turned out relevant,
 see Theorem \ref{t:kpcqpq} and its subsequent remark,
 keeping also Theorem \ref{t:chi-ap} in mind:

\bpm
Investigate the role of edge connectivity number in the context
 of\/ $g(n,G|H)$ for graphs\/ $H$ which admit a bipartite anti-packer
 inside\/ $G$.
\epm

Consider the largest possible value of $d$ in a feasible pair, $d=2$.
In the non-degenerate case the \tur\ function depends only
 on the chromatic number, therefore it
 may even hold that for $d=2$ all feasible $(c,d)$ can be written
 in the form $(\frac{r-2}{2r-2},2)$ for some integer $r\geq 3$.
This question is open so far in general, nevertheless
 Theorem \ref{t:chi-ap} specifies a large class of instances for
 which it is proved that other pairs $(c,d)$ with $d=2$ cannot exist.

The other extreme, $d=0$ represents $\phi(n,G|\cF)=O(1)$.
In this case we can prove that (\ref{eq:limit}) is satisfied.

\btm
If\/ $\phi(n,G|\cF)=O(1)$, then there is a constant\/ $c$ and a
 threshold\/ $n_0=n_0(G,\cF)$ such that\/ $\phi(n,G|\cF)=c$
 for all\/ $n>n_0$.
\etm

\bpf
Let $p:=|G|$.
Assuming $\phi(n,G|\cF)=O(1)$, there exists a universal
 upper bound say $k$ on $\phi(n,G|\cF)$.
Then also the isomorphism types of edge colorings of $K_p$ with
 at most $k$ distinguished colors (not allowing the renumbering
 of colors) is bounded, say it is at most $m(k)$.

We let $n'=m(k)\cdot p + 2k$.
With this choice, the following properties are valid in any
 edge coloring $\psi$ of $K_n$ for any $n>n'$:
 \begin{enumerate} {\it 
  \item {\rm If more than $k$ colors are used, then there exists a
   vertex $w$ such that $K_n-w$ still contains at least $k$ colors.}
  \item {\rm If at most $k$ colors are used, then there exists a
   vertex $w$ such that $K_n-w$ still contains the same set of
   color patterns induced by its copies of $K_p$ as $K_n$ does.
    And, in particular, all colors of $\psi$ are present
     also in $K_n-w$.}}
 \end{enumerate}
\emph{Proof of 1.:} 
Assume that the number of colors used is $k+t$ for some $t\geq 0$.
If a color class is a single edge or a star, mark the two ends
 of that color or the center of the corresponding star with
 the color.
There are at most $2k+2t$ marks on the total of $n$ vertices,
 thus a vertex $w$ marked with at most
 $$
   \floor{\frac{2k+2t}{n}} \leq \floor{\frac{2k+2t}{2k+p}} \leq 
    1 + \floor{\frac{t-p/2}{k+1}} \leq t
 $$
  exists.
Removing this $w$ does not decrease the number of colors to
 less than~$k$.

\msk

\nin
\emph{Proof of 2.:}
If $k'$ is the number of colors, then in the current case we have
 $k'\leq k$ and $m(k')\leq m(k)$.
For each of the $m(k')$ possible coloring patterns $\pi$ of $K_p$,
 select a $p$-element subset $X_\pi$ in which $\psi$ induces $\pi$,
 if $\pi$ occurs in $K_n$.
We set $X:= \bigcup_\pi X_\pi$.
By the choice of $n$ we have $n>|X|$.
Pick any vertex $w\in V(K_n)\smin X$.
By the choice of $X$, all
 patterns incident with $w$ have their copies in $K_w$.

\msk

If $\phi=f$, let $\psi$ be any edge coloring of $K_n$ using
 $k'\geq f(n-1,G|\cF)$ colors for $K_n$.
If $k'\geq k$, a copy of $G$ required for $f$ exists by the
 choice of $k$.
If $f(n-1,G|\cF) \leq k' < k$, we remove $w$ as guaranteed by {\it 2.}
Since at least $f(n-1,G|\cF)$ colors remain, a required copy of
 $G$ is present; and as this property holds also for
 $k'=f(n-1,G|\cF)$, we obtain $f(n,G|\cF)\leq f(n-1,G|\cF)$.

If $\phi=g$, consider an edge coloring $\psi$ of $K_n$ with
 $k':=g(n,G|\cF)$ colors, such that $\psi$ realizes the
 requirements of $g(n,G|\cF)$.
Then of course $k'\leq k$, by the choice of $k$.
Applying {\it 2.}, we obtain a coloring of $K_n-w\cong K_{n-1}$,
 still with $k'$ colors, in which all copies of $G$ contain a
 monochromatic copy of some $F\in\cF$.
This implies $g(n-1,G|\cF)\geq g(n,G|\cF)$.

All the above imply that $(k_n)_{n\geq 1}$ with $k_n:=\phi(n,G|\cF)$
 is an ultimately non-increasing sequence of nonnegative integers.
Thus, it is a constant for all $n>n_0$, and (\ref{eq:limit}) is valid.
\epf

\bpm
Characterize the pairs\/ $(G,\cF)$ for which\/
 $\phi(n,G|\cF) = O(1)$ as\/ $n\to\infty$.
\epm

\bpm
Given an integer\/ $c>0$, characterize the pairs\/ $(G,\cF)$
 for which\/ $\phi(n,G|\cF) = c$ for all\/ $n$ sufficiently large.
In particular, is\/ $(k,0)$ a feasible pair for every positive
 integer\/ $k$\,?
\epm

Perhaps the second part of this problem is not too hard if $k$ is
 very small.

Let us formulate a more explicit version of Problem \ref{pr:fi-diff}.

\bpm
Are the feasible pairs the same for\/ $f$ and for\/ $g$?
\epm

It seems plausible to guess that the answer is affirmative.
A partial evidence
 is the following fact for $d>0$, by which this problem
 seems to be much easier than the complete characterization of
 feasible pairs.

\bpn   \label{p:compl-c,d}
Assume that\/ $\phi(n,G|\cF)\to\infty$ as $n$ gets large, and\/
 $(G,\cF)$ determines the feasible pair\/ $(c,d)$.
If\/ $K_2\in \cF$ and\/ $G\notin \cF$, or if\/ $K_2\notin \cF$ and\/
 $G\in \cF$, then\/ $(c,d)$ is a feasible pair also for\/
 $\{f,g\}\smin\{\phi\}$.
\epn

\bpf
This follows by Proposition \ref{p:compl}, taking the family
 complementary to $\cF$.
\epf

Let us recall Proposition \ref{p:monoton}, which implies that
 if $\cF\subset \cF'\subset \cF''$ are any three nested families,
 such that $(G,\cF)$ and $(G,\cF'')$ determine the same pair
 $(c,d)$, then also $(G,\cF')$ determines the same pair.
This fact motivates the following natural problem.

\bpm
Given a graph\/ $G$ and the real numbers\/ $c,d$ with\/ $c>0$ and\/
 $0\leq d\leq 2$, find the inclusionwise minimal and maximal
 families\/ $\cF$ such that\/ $G$ and\/ $\cF$ determine\/ $(c,d)$.
\epm

Although not minimal or maximal, nevertheless many ``small'' and
 ``large'' families can be observed in the three parts of our trilogy.
In this way the interested reader can state various ``sandwich
 theorems'', mostly in the quadratic scenario.
Here we give only one example.

\bex
Let\/ $G=K_p$, $p\geq 5$, and let\/ $\phi=f$.
Search for families\/ $\cF$ such that\/ $(G,\cF)$ generate\/
 $(c,d)=(\frac{r-2}{2r-2},2)$ with\/ $r=\ceil{p/2}$.
It is equivalent to the asymptotics\/
 $f(n,G|\cF) = (1+o(1))\,\ex(n,K_{\ceil{p/2}})$.

The extremely dense family\/ $\nem{P_3}$ obtained by the omission
 of the single graph\/ $P_3\cong K_{1,2}$ satisfies the conditions,
 undoubtedly it is a very transparent maximal family.

From the other side, there are sparse families with the same\/
 $f(n,G|\cF) = (1+o(1))\,\ex(n,K_{\ceil{p/2}})$.
Some examples are:\/
 $\fsq$ (the graphs with square number of edges),\/
 $\cF_{1,k}$ (the graphs whose number of edges is congruent to\/
  $1$ {\rm modulo} $k$,
 or\/ $\{mK_2\mid m\geq 1\}$ (matchings), the latter being the
  defining family of the local \arr\ numbers\/ $\lar(n,G)$.
These families seem to be rather unrelated, so perhaps it is
 rather hard to characterize the minimal classes generating
 this partiular\/ $(c,d)$.

Recall further that\/ $\nem{K_{1,t}}$ leads to the analogous function\/
 $\ex(n,K_{\ceil{p/t}})$, which means among other things that in case of\/
  $p/t=p'/t'$ we arrive at the same asymptotics for\/
  $f(n,K_p|\nem{K_{1,t}})$ and\/ $f(n,K_{p'}|\nem{K_{1,t'}})$.
This is an additional indication that the problem of minimality
 for\/ $(c,d)$ is really hard.
\eex

Note that this ``sandwich'' phenomenon is brand new in \arr\
 theory.
This possibility emerged only in the context of the functions $f$
 and $g$, because the classical concept of $\AR(n,G)$ is defined
 purely in terms of $\cF=\{K_2\}$.

\subsection{Miscellaneous open problems}

Concerning $f(n,K_p|\nem{P_3})$ the smallest possible instance is
 $$
   f(n,K_3|\nem{P_3}) = n \,.
 $$
The upper bound is implied by
 $f(n,K_3|\nem{P_3})\leq \AR(n,K_3)=n$, whereas the lower bound
 is implied by the fact that any three vertices in the \Lex\
 pattern (using $n-1$ colors) induce
 a monochromatic $P_3$ as a color class, hence not feasible
 for $\nem{P_3}$.
However, already the second ssmallest case is far from being solved.

\bpm
Prove or disprove that\/ $f(n,K_4|\nem{P_3})$ is linear in\/ $n$.
\epm

So far our best lower bound on $f(n,K_4|\nem{P_3})$ is $2n-2$,
 obtained by a rainbow star $K_{1,n-1}$ together with the
 \Lex\ pattern in the $K_{n-1}$ subgraph induced by its leaves.
Our currently best upper bound grows with $O(n^{3/2})$,
 derived from the presence of a clean $K_{2,2}$ via Theorem \ref{t:CCCL}.
In fact the local \arr\ function $\lar(n,K_4) = \Theta(n^{3/2})$
 is an upper bound on $f(n,K_4|\nem{P_3})$, but the
 lower-bound construction on $\lar$ is not applicable here.
We note that $2n-2$ is the right value for $n=4$, because
 a monochromatic $P_3\subset K_4$ and four private colors
 on the edges of $K_4 - P_3$ yield a 5-coloring
 that does not satisfy the requirements on $f$.

\bpm
Determine\/ $\displaystyle \lim_{n\to\infty}
  \frac{f(n,K_{2p}|\upp{pK_2})}{n}$ and\/ $\displaystyle
   \lim_{n\to\infty} \frac{g(n,K_{2p}|\upp{pK_2})}{n}$.
\epm

\bpm
Determine\/ $g(n,K_4|2K_2)$, and more generally\/
 $g(n,K_{2p}|pK_2)$ and beyond,\/ $g(n,K_p|qK_2)$ where\/ $p>2q$.
\epm

\bpm
Determine the asymptotics of\/ $g(n,K_p|\upp{P_q})$ for\/ $p\leq 2q-2$.
\epm

\bpm
Determine the asymptotics of\/ $g(n,K_p|\upp{C_q})$ for\/ $p\leq 2q-2$.
In particular, settle the cases\/ $(p,q)=(5,4)$ and\/ $(p,q)=(6,4)$.
\epm

\bpm
Determine the asymptotics of\/ $g(n,K_p|\upp{K_{1,p-1}})$, and
 more generally of\/ $g(n,K_p|\upp{K_{1,q}})$ for\/ $p\leq 2q-2$.
\epm

Although more ambitious, it may even be possible to determine the
 exact values of some of $g(n,K_p|\upp{P_q)}$, $g(n,K_p|\upp{C_q)}$,
 or $g(n,K_p|\upp{K_{1,p-1}})$, for certain ranges of $p$
 with respect to $q$, or at least for some pairs $p,q$ beyond
 those proved in Section \ref{s:pierc}.

\bpm
Investigate the behavior of\/ $g(n,K_p|\upp{T})$ where\/ $T$ is a
 forest other than star, path, or matching.
\epm

\bpm
Investigate the behavior of\/ $g(n,K_p|\upp{F})$ where\/ $F$ is a
 graph of small order.
\epm

\bpm
Investigate\/ $f(n,K_p|\fs{S})$ for the sequences\/ $S = S(r,k) =
 \{1\} \cup \{ m \mid m\equiv r \ (\mathrm{mod} \ k)\}$, where\/
 $r\neq 1$.
\epm

\bpm
Determine\/ $\AR(n,C_{2k+1}\cup K_2)$, and more generally\/
 $\AR(n,C_{2k+1}\cup mK_2)$.
\epm

Concerning this problem it is worth recalling the important role
 of $C_{2k+1}\cup K_2$ in the determination of
  $g(n,p,q) = g(n,K_p|\upp{K_q})$ for $p=2q-2$.
In fact our estimates on $\AR(n,\{C_{2k+1}\cup K_2,C_{2k+2}\})$
 are fairly tight.
However, the analogous question for the other parity remains open.

\bpm
Determine\/ $\AR(n,\{C_{2k}\cup K_2,C_{2k+1}\})$.
\epm

\bpm
Determine\/ $\AR(n,\cH_{p,d})$ for\/ $p>d>1$, where\/ $\cH_{p,d}$ is
 the family of graphs of order\/ $p$ and minimum degree\/ $d$.
\epm

A motivation for this problem is that the members of $\cH_{p,d}$
 are blockers of $K_{1,q}$ with $q=p-d+1$ whenever $1<q<p$.

\end{document}